\documentclass[a4paper,11pt]{amsart}

\usepackage{amsmath,amssymb,amsfonts,latexsym}
\usepackage[english]{babel}
\usepackage[dvipsnames]{xcolor}
\usepackage{multirow}
\usepackage[left=2.5cm,right=2.5cm,top=3.cm,bottom=3.cm]{geometry}
\usepackage{tikz}
\usepackage{pgfplots}
\usepackage{enumitem}
\usepackage{textcomp}
\usepackage{colortbl}
\usepackage{slashbox}
\usepackage{setspace}
\usetikzlibrary{pgfplots.statistics}
 
\usepackage{arydshln}
\usetikzlibrary{matrix,shapes,arrows,positioning}
\definecolor{myblue}{RGB}{209,221,243}
\definecolor{myblu}{RGB}{110,148,216}
\definecolor{myre}{RGB}{255,188,188}
\definecolor{li}{gray}{0.97}
\definecolor{lt}{gray}{0.7}
\usepackage{tabu}
\usepackage{pifont}
\usepackage{hhline}
\usepackage{algorithm}
\usepackage[noend]{algpseudocode}
\makeatletter
\def\BState{\State\hskip-\ALG@thistlm}
\usetikzlibrary{shadows}
\usepackage{booktabs}

\usetikzlibrary{patterns}
\newcommand{\capt }[2]{\textbf{\vspace{3mm} \\ } T{\footnotesize ABLE} #1. #2}

\algrenewcommand\algorithmicindent{1.4em}%

\begin{document}

\title[Parallel MOEA with Criterion-based Selection applied to the KP]{\textbf{\vspace{-3mm} \\}A Parallel MOEA with Criterion-based Selection Applied to the Knapsack Problem  }

\author[KANTOUR Nedjmeddine, BOUROUBI Sadek and Chaabane Djamel]{ KANTOUR Nedjmeddine$^1$, BOUROUBI Sadek$^2$ and Chaabane Djamel$^3$ \vspace{2mm} \\ \MakeLowercase{\texttt{nkantour@usthb.dz}$^1$, \texttt{sbouroubi@usthb.dz}$^2$, \texttt{dchaabane@usthb.dz}$^3$}  \vspace{3mm}\\\emph{$^{1,2\,}$LIFORCE L\MakeLowercase{aboratory}, $^3\,$L\MakeLowercase{aboratory} AMCD-RO \vspace{1mm}\\
USTHB, F\MakeLowercase{aculty of }M\MakeLowercase{athematics}, D\MakeLowercase{epartment of} O\MakeLowercase{perations }R\MakeLowercase{esearch},\vspace{1mm} \\ P.B. 32 E\MakeLowercase{l-Alia}, 16111, B\MakeLowercase{ab }E\MakeLowercase{zzouar}, A\MakeLowercase{lgiers}, A\MakeLowercase{lgeria.}}
}

%
%



\begin{abstract}
 In this paper, we propose a parallel multiobjective evolutionary algorithm called Parallel Criterion-based Partitioning MOEA (PCPMOEA), with an application to the Mutliobjective Knapsack Problem (MOKP).  The suggested search strategy is based on a periodic partitioning of potentially efficient solutions, which are distributed to multiple multiobjective evolutionary algorithms (MOEAs). Each MOEA is dedicated to a sole objective, in which it combines both criterion-based and dominance-based approaches. The suggested algorithm addresses two main sub-objectives: minimizing the distance between the current non-dominated solutions and the ideal point, and ensuring the spread of the potentially efficient solutions. Experimental results are included, where we assess the performance of the suggested algorithm against the above mentioned sub-objectives, compared with  state-of-the-art results using well-known multi-objective metaheuristics. \\


\vspace{2mm}

\vspace{-2em}

\end{abstract}



\maketitle

\section{\textbf{Introduction}}
Multiobjective Optimization Problems (MOPs) involves several conflicting criteria, where any solution is qualified as optimal if it belongs to criteria trade-offs set called the Pareto front \cite{30}. In other words, each Pareto optimal solution reaches a good tradeoff among these conflicting objectives: one objective cannot be improved without worsening
any other objective \cite{28}. The difficulty of multiobjective optimization lies in the absence of a total order relation which links all the solutions of the problem. In terms of evolutionary algorithms, this lack appears in the difficulty of conceiving a selection operation that assigns for each individual a probability of selection proportionally to the performance of that individual. Another drawback is the premature loss of diversity; hence the need to conceive techniques for maintaining diversity within the population. Against such problems one have to conceive algorithms that satisfies the following criteria \cite{29}: the algorithm must converge towards the true Pareto front in a reasonable time, and it must also result diverse solutions on the front in order to have a good representative sample instead of focusing on an area of the objective space.
 A number of approaches (exact and heuristic) have been considered for solving MOPs, mainly, metaheuristics which have proven their ability to give good approximations. Among these methods, one of the most popular resorts for solving MOPs are Multiobjective Evolutionary Algorithms (MOEAs). The very early works on MOEAs were initiated by D. Schaffer (VEGA, \cite{8}), who is considered to be the first to design a MOEA during the 1980s. Since than, many successfull MOEAs has been proposed, see for exapmle: PAES \cite{41}, NSGAII \cite{6}, SPEA2 \cite{7}, IBEA \cite{44},  MOEA/D \cite{37}, etc. The reader may find further details about early studies of MOEAs in \cite{30}. The employment of MOEAs in Multiobjective Optimization (MOO) is often justified by their population based character, managing to achieve high quality approximations of Pareto optimal solutions.

On the other hand, For most MOPs, executing the generational cycle of standard MOEAs on large instances of the problem and/or on large sizes of population requires considerable resources in terms of computational time and memory \cite{1}. Therefore, a variety of design and implementation difficulties are studied to construct more effective MOEAs. Overcoming these difficulties usually involve defining new operators, hybrid algorithms, and parallel models, etc. Hence, the parallelization of MOEAs emerges naturally when dealing with computationally expensive algorithms. However, parallelization of MOEAs aims not only to reduce computational time, but also to improve the quality of the approximated Pareto fronts, and increases the robustness of MOEAs against MOPs in both real life and theoretical research fields \cite{1}. In this work, we employed the master/worker model, handling multiple MOEAs with different populations, and a separated "main" monitoring algorithm. Many pMOEAs of different models and implementations can be found in literature, for literature review, see \cite{2}, \cite{34}, \cite{35}. However, in the following part we will briefly present some of the existing pMOEA using master/worker model and/or multiple search algorithms (MOEAs), partitioning either the search space or the objective space.

\subsection*{Related work}  
An early work on master/worker pMOEA named the Parallel Multi-Objective Genetic Algorithm
(PMOGA) applied to eigenstructure assignment problems \cite{31}. This algorithm launches multiple MOEAs using identical initial population and different decision-making logic (fitness assignment and selection operator), where each one is executed on a worker processor. The master uses the populations produced by workers to form final population \cite{2}. An other master/worker pMOEA is applied to solve scheduling problems \cite{32}, where two versions were implemented: heterogeneous and homogeneous populations. The first uses multiple subpopulations distributed for each objective function, while the heterogeneous subpopulations are Pareto-oriented. Furthermore, both versions utilize unidirectional migration flow by sending individuals to a separate main population~\cite{2}. We mention also the master/worker Parallel Single Front Genetic Algorithm (PSFGA) applied to several the benchmark functions \cite{32}. This algorithm performs periodically the following tasks: it sorts the population according to the values of the objective functions, then it partitions it into subpopulations and send them to different processors, where a sequential multiobjective genetic algorithm (SFGA) is applied to each subpopulation \cite{32}. In \cite{22}, the authors propose an extended version of IBMOLS to a cooperative model with an application to the multiobjective  Knapsack Problem.  The suggested approaches uses a multipopulation based cooperative framework W${\footnotesize \epsilon}$-CMOLS, in which different subpopulations are executed independently with different configurations of IBMOLS. Furthermore, each subpopulation is focused on a specific part of the search space using a weighted versions of the epsilon indicator \cite{22}. Finally, a hybrid swarm intelligence algorithm called Multi-objective Firefly Algorithm with Particle Swarm Optimization (MOFPA) applied to the multiobjective Knapsack Problem \cite{5}. This algorithm combines the Firefly Algorithm (FA) and Particle Swarm Optimization (PSO), and it employs a $\tanh$ transfer function for the discretization of these algorithms; furthermore, it uses the epsilon dominance relation for managing the size of the external archive (an external non-dominated population) \cite{5}. 

In this paper, we propose a parallel multiobjective evolutionary algorithm designed in a master/worker model, we call it Parallel Criterion-based Partitioning MOEA (PCPMOEA). The proposed algorithm uses multiple MOEAs with a criterion-based selection operation, each one is dedicated to a sole objective function, while these algorithms are monitored by a master entity that periodically collects and partitions the current Pareto solutions according to their distribution in the objective space. The parts are used to update the current subpopulations of MOEAs. Experimental results are provided, where we used several benchmark instances of the multiobjective Knapsack Probkem, and compared PCPMOEA with some effective multiobjective evolutionary algorithms. This paper is organized as follows: in the next section, we recall some basic concepts of MOO, and the multiobjective Knapsack Problem. In section 3, we present a brief description of parallel MOEAs and its taxonomy. In section 4, we present the suggested parallel approach, with detailed functioning description. In section 5, we discuss the experimental results obtained from the suggested algorithm. Finally, in Section 5, we end the paper with the conclusion and our perspectives.\vspace{-3mm}
\section{\textbf{Background}}
\subsection*{Multiobjective Problems} The goal in Multiobjective Problems (MOPs) is to optimize simultaneously $k$ objective functions. More formally, MOPs can be stated as follows \cite{3}:
$$ (MOP)
\left\{
\begin{array}{l}
  "max"\, Z(x) = (Z^1(x), Z^2(x), \dots, Z^k(x)), \\
  \; s.t. \quad x \in \Omega. 
\end{array}
\right.
$$
where $\Omega$ is the decision space, and $x \in \Omega$ is a decision vector (feasible solution). The vector $Z(x)$ consists of $k$ objective functions.
$$Z^i : \Omega \rightarrow \mathbb{D}_i,\; i = 1, \dots, k,$$ 
where, $\displaystyle \prod_{i=1}^n\mathbb{D}_i $ is the feasible objective space.\\

Very often, the objectives in MOPs contradict each other, there exists no point that maximizes all the objectives simultaneously. One has to balance them. Moreover, a tradeoffs solutions set can be defined in terms of Pareto optimality using the dominance relation.
\subsection*{Pareto optimality}
Since the aim in MOPs is to find good compromises rather than a single solution as in single-objective optimization problems, we present the dominance relation, allowing to define optimality for MOPs. A solution $x$ is said to dominate another solution $x^{'}$, denoted as $x \succ x^{'}$, if and only if, $\forall i \in \{1, \dots , k\}, Z^i(x) \geq Z^i(x^{'})$ and $Z(x) \not = Z(x^{'})$. A feasible solution $x^* \in \Omega$ is called a Pareto optimal solution, if and only if, $\not\exists y \in \Omega$ such that $Z(y)\succ Z(x^{*})$. The set of all Pareto optimal solutions is called the Pareto-optimal set (PS); moreover:
$$PS = \{ x \in \Omega | \not\exists y \in \Omega, Z(y)\succ Z(x^{*})\}. $$
 The evaluation of PS in the objective space is called the Pareto front (PF):
$$PF = \{Z(x)| x\in PS \}. $$

From a MOEA point of view, Pareto-optimal solutions can be seen as those solutions within the genotype search space, whose corresponding phenotype objective vector components cannot be all simultaneously improved \cite{2}.
\subsection*{Ideal vector}
The ideal vector contains the optimum for each separately considered objective function, all constituing the same vector in the objective space (usually $\mathbb{R}^n$) \cite{2}.
In other words, let $x^{(i)}_0 \in \Omega$ be a given vector of decision variables which maximizes the $i^{th}$ objective function $Z^i$. Meaning that the vector $x^{(i)}_0$ verifies the following equality \cite{2}:
$$ Z^i(x^{(i)}_0) = \underset{x \in \Omega}{\mbox{max}} \, Z^i(x).$$

We denote by $Z^{i}_0$ the optimum evaluation of the $i^{th}$ objective function $Z^i(x^{(i)}_0)$, and by the vector $  Z_0 = (Z^{1}_0, Z^{2}_0, \dots, Z^{k}_0) $ the ideal  vector for a MOP \cite{2}. The ideal vector is usually used in some methods as a reference point. 

\subsection*{\textbf{Multiobjective multidimensional Knapsack Problem}}
 The Multiobjective multidimensional Knapsack Problem (MOMKP) is a widely studied combinatorial optimization problem. Furthermore, the MOMKP is the multidimensional version of the multiobjective Knapsack Problem (MOKP) \cite{3}, which is proven to be $\mathcal{NP}-$hard \cite{15}, besides, it is known for the fact that the size of the Pareto optimal solutions set can grow exponentially with the number of items in the knapsack \cite{20}. Despite its simple Formulation, MOKP can be applied for modeling many real problems such as budget allocation and resource allocation~\cite{5}. Therefore, its resolution has both theoretical and practical character that draws researchers attention. MOMKP is a particular case of multiobjective linear integer programming (MOILP) \cite{3}. Mathematically, MOMKP can be stated as follows: given $n$ items having $p$ characteristics $w^i_{j}, j\in \{1, \dots , n\}$ (weight, volume, etc.), and $k$ profits $c^{i}_{j}, j\in \{1, \dots ,n\}$, we want to select items as to maximize the $k$ profits, while not exceeding the $p$ knapsack capacities $W_i$ \cite{3}.
$$(MOMKP)
\left\{
\begin{array}{l}
\vspace{-0.5cm}\\
 \vspace{-5mm} \\
\displaystyle "max\," \ Z^i(x)=\sum_{j=1}^{n} c^{i}_{j}x_{j},\quad i\in \{1, \dots ,k\}\vspace{-1mm}\\
\begin{array}{ll}
s.t. & \displaystyle \sum_{j=1}^{n} w^i_{j}\,x_{j} \leq W_i ,\quad i\in \{1, \dots ,p\}\vspace{0.4em}\\
 & {\displaystyle x_{j} \in \{0,1\} ,\: \forall j\in \{1, \dots , n\}, }\vspace{1.3mm}
\end{array}
\end{array}
\right.\vspace{3mm}
 $$ 
 where $n$ is the number of item, $x_j$ denotes the $j^{th}$ binary decision variable, $k$ represents the number of objectives, $Z^i$ stands for the $i^{th}$ objective function.
  
 The multiobjective Knapsack Problem (MOKP) and its variants have been subject to many studies, addressing new approaches and analysis of either exact or heuristic resolution methods. Namely, Zitzler and Thiele \cite{16} pioneered the work on the multidimensional variant of the MOPK using evolutionary algorithms, where they introduced  also a set of MOMKP instances that are widely used afterwards. Many other research works are dedicated to different variants of the MOKP (see for example \cite{17}, \cite{18}, \cite{19}, \cite{20}, \cite{24}, \cite{27}). However, the existing heuristic approaches are not restricted to evolutionary algorithms, one may find for instance: indictor based ACO \cite{21}, parallel multi-population algorithm using local search \cite{22}, iterated local search based on quality indicators \cite{23}, tabu search \cite{25}, and swarm intelligence \cite{5}.  The interested reader may find a detailed discussion of resolution approaches of the MOKP and its variants in \cite{3}. Due to the popularity and simplicity of the MOKP, and the fact that MOEAs are known to perform at their best against problems with string structures; we have chosen the MOMKP as test problem for the assessment of the suggested approach. Furthermore, we used the MOMKP known benchmark instances itroduced by Zitzler and al. in \cite{16} (also referred to ZMPK instances in \cite{3}). \vspace{-2mm}
\section{\textbf{Parallel multiobjective evolutionary algorithms (pMOEAs)}}
In this section we present a brief sweep over the common theoretical aspects of parallelization models and taxonomy that are used in (MO)EA. MOEAs are the most commonly considered metaheuristics for parallelization. That is when dealing with populations of individuals, parallelism arises naturally (i.e., each individual can be considered as an independent unit) \cite{13}. As a natural extension, the parallelization of MOEAs are derived from models designed for single-objective optimization \cite{12}: master/worker, island, and diffusion models. Furthermore, there exists many criteria that provides the theoretical description of a given parallelization.  In a recent publication, T. El-Ghazali has given a taxonomy for pMOEAs, where it describes a unified view over pMOEAs \cite{1}.

 However, in this work we adopt a master/worker or global parallelization paradigm. Where a processor serves as master, which usually preforms tasks that require a global outlook over the target problem (search space, objective functions, etc.), as it is the case for selection operator in MOPs. Furthermore, the master supervises multiple processors called workers, by distributing workload (i.e., tasks, subproblems) to the workers to process (i.e., execute, solve). Moreover, for the design of pMOEA there are also three major hierarchical models that can be identified \cite{1}: algorithmic-level, iteration-level, solution-level. This classification is done according to their dependency on the target MOPs, behavior, granularity, and aim of parallelization. Note that each paradigm can be implemented in either a synchronized or asynchronized way \cite{14}.  A further particularity in designing a pMOEA is to specify whether the whole Pareto front is distributed among each parallel search algorithm, or it is a centralized element of the parallel scheme. This leads also to another classification criterion of pMOEAs, where these two strategies are called: distributed and centralized approaches \cite{1}, \cite{13}.

Accordingly, the suggested parallel algorithm is designed in a master/worker paradigm with full independence from the target multiobjective problem. Besides it works in an algorithmic-level parallelization, where we use multiple asynchronous search algorithms working on sub-populations in parallel (intra-algorithm parallelization); each one is assigned a region of the objective space to work on. Furthermore, the suggested search strategy uses a centralized Pareto front handled by the master. This Pareto front is built by combining multiple distributed Pareto fronts where "local" non-dominated solutions are archived. Then, this latter is advisedly partitioned and distributed again to workers/search algorithms as to proceed "local" search. 
\section{\textbf{Parallel Criterion-based Partitioning MOEA (PCPMOEA) }} 
In this section, we describe the general functioning of the suggested parallel scheme, and the search strategy. The idea is to launch multiple asynchronous MOEAs with different populations, where each MOEA is dedicated to a given criterion. While the search entities (MOEAs) are supervised by a global processing element, in order to adjust and redirect the search process in real time. In summary, The suggested pMOEA can be classified as a cooperative algorithmic level parallel model designed in a master/worker paradigm, handling:\vspace{1mm}
\begin{itemize}[leftmargin=0.2in]
\item Workers (search entities) consisting of criterion-based MOEAs: multiple MOEAs with criterion-based selection.\vspace{1mm}
\item Master entity where its main mission is to update the current population of each search entities (workers). That is, by preforming a global selection among elite individuals produced by the workers, and then advisedly partitions and redistributes the selected individuals to search entities.\vspace{2mm}
\end{itemize}
\begin{figure}[H]
\begin{tikzpicture}[transform shape, scale =1.]
\hspace{8mm}
\tikzstyle{debutfin}=[pattern=north west lines, pattern color=gray!40,rectangle,rounded corners=.9pt, minimum width=1.6cm,
                        minimum height = 0.7cm,draw]
\tikzstyle{debutfi}=[fill=gray!20,rectangle,rounded corners=.9pt, minimum width=1.6cm,
                        minimum height = 0.7cm,draw]
\node[debutfin] (a) at (-1.5,0.6) {$Master$};
\node[debutfi] (b) at (-5,3) {$MOEA_1$};
\node[debutfi] (c) at (-2.6,3) {$MOEA_2$};
\node (d) at (-.7,3) {{\Huge \dots}};
\node[debutfi] (e) at (1.4,3) {$MOEA_k$};

\draw[-, thick, black] (-1.8,.95) -- (b.south);
\draw[-, thick, black] (a.north) -- (c.south);
\draw[-, thick, black] (-1.2,0.95) -- (e.south);

\draw[<->,very thick, black!60] (2.4,.9) -- (2.4,2.5) node[midway,sloped,right,rotate=270] {\color{black}\begin{tabular}{c}
Information\\ exchange
\end{tabular}
};
\end{tikzpicture}
\caption{Schematic of the adopted parallel model.}
\label{fig:pscheme}
\end{figure}
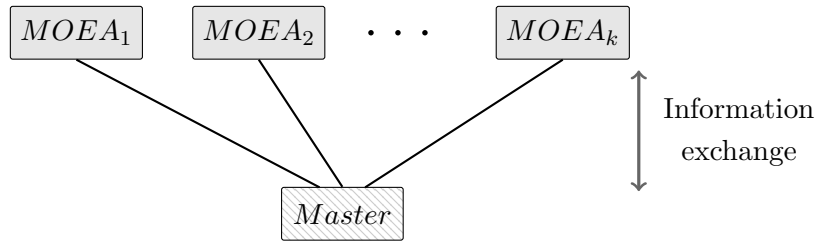

In other words, each $MOEA_i, i=1, \dots, k$, represents a criterion based search entity, dedicated to the $i^{th}$ objective function. The master node accomplishes (periodically) the following tasks: performs global selection among the potentially efficient solutions produced by the search entities, and then updates the current population in each MOEA (see Figure \ref{fig:pscheme}). In the rest of this section, we will present a detailed description of the suggested pMOEA.

\subsection{\textbf{$MOEA{\Large_{i}}$ (criterion-based MOEA)}} We want to design a MOEA that aims to minimize the distance between the current Pareto front and the ideal point. That is by focusing each MOEA on a single objective, while retaining its multiobjective character. The suggested MOEA is characterized by a secondary population and its criterion-based selection mechanism. Many MOEAs involve an additional population in the search process, as to store potentially efficient solutions found through the generational process. This approach has proven to be effective in finding good approximations of the optimal Pareto front (see for example: PAES \cite{41}, SPEA \cite{16}, SPEA2 \cite{7}, MOMGA \cite{42}, MOMGA-II  \cite{43}). This secondary population is often called archive or external archive \cite{2}. In this algorithm, like most of MOEAs, we follow the pattern of initializing a population of individuals, then executing a generational loop with evolutionary operators, fitness assignment, ranking individuals, and then storing non-dominated solutions in an archive. Furthermore, the archive size is not limited, while the number of dominated individuals is fixed beforehand. However, the archive is updated periodically by the master entity as to keep it enclosed and dedicated to its allocated search space. The following algorithm resumes the adopted MOEA, we call it $MOEA_i$ referring to the $i^{th}$ objective function $Z^i$:

\begin{algorithm}
{\fontsize{11pt}{11pt}\selectfont
\caption{$MOEA{\Large_{i}}$ (criterion-based: $i^{th}$ objective)}\label{alg:HI}
\begin{algorithmic}[l]
\State \vspace{0.2em}
\State \hspace{-0.5cm}\textbf{Begin} \vspace{.15cm}
\State Generate the initial population $P_{i}$;
\vspace{0.15cm}
\Repeat \vspace{0.15cm}
\State Evaluate the individuals in $P_{i}$;\vspace{0.15cm}
\State Apply selection operator (see Figure \ref{fig: eai});\vspace{0.15cm}
        \State Select parents form $P_{i}$ and perform evolutionary operations;\vspace{0.15cm}
        
        \If{(Migration condition is satisfied)\vspace{0.15cm}}
            \State Send the current archived solutions to master;\vspace{0.15cm}
            \State Update archive (see Section \ref{master});\vspace{0.15cm}
            \State \hspace{-7.5mm} \textbf{endif}\vspace{0.15cm}
        \EndIf
        \vspace{0.15cm}           
       
      \Until{(Termination condition met)\vspace{0.15cm}}\vspace{0.1cm}
\State \hspace{-0.5cm}\textbf{End} \vspace{0.4em}
\end{algorithmic}
}
\end{algorithm}\vspace{-3mm}
\subsubsection*{\textbf{Selection mechanism}}
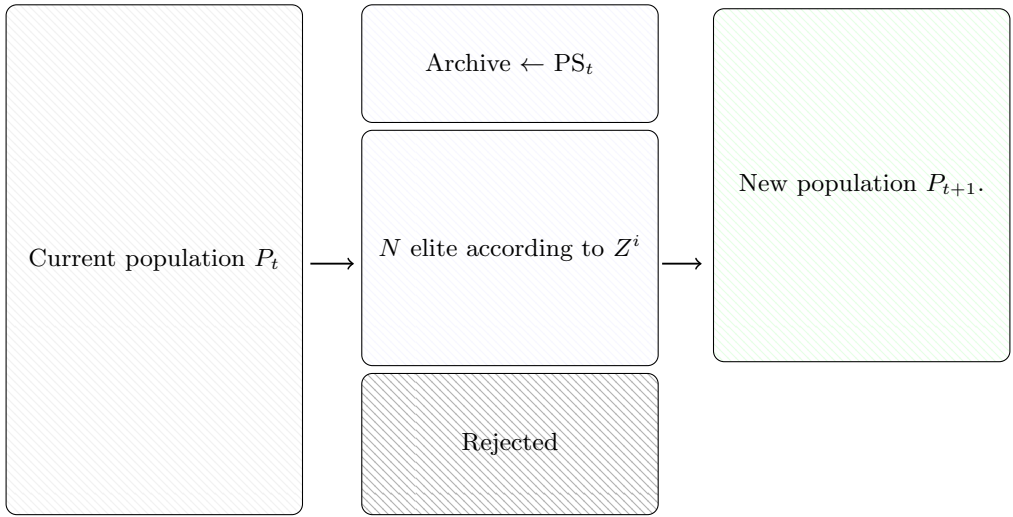
\begin{figure}[h!]
\tikzset{
   brace/.style={
     decoration={brace, mirror},
     decorate
   }
}

\begin{center}
\begin{tikzpicture}[transform shape, scale = 1.04]

\draw[->, line width=.25mm] (2.85,9.2) -- (3.4,9.2);
\draw[->, line width=.25mm] (7.3,9.2) -- (7.8,9.2);

\draw [pattern=north west lines, pattern color=gray!10, rounded corners=4pt] (-1,6) rectangle++ (3.75,6.5) node[pos=.5] (a)  {{\footnotesize Current population $P_t$}};

\draw [pattern=north west lines, pattern color=blue!3,rounded corners=4pt] (3.5,11) rectangle++ (3.75,1.5) node[pos=.5] (b)  {{\footnotesize Archive $\leftarrow$ PS$_t$}};

\draw [pattern=north west lines, pattern color=blue!3, rounded corners=4pt] (3.5,7.9) rectangle++ (3.75,3) node[pos=.5] (c) {{\footnotesize $N$ elite according to $Z^i$}};

\draw [pattern=north west lines, pattern color=gray!60,rounded corners=4pt] (3.5,6) rectangle++ (3.75,1.8) node[pos=.5]  {{\footnotesize Rejected}};

\draw [pattern=north west lines, pattern color=green!10, rounded corners=4pt] (7.95,7.95) rectangle++ (3.75,4.5) node[pos=.5]  {{\footnotesize New population $P_{t+1}$.}};
\end{tikzpicture}\vspace{-2mm}
\end{center}
\caption{Schematic of the selection mechanism (iteration $t$).}
\label{fig: eai}
\end{figure}

The most commonly employed selection operators are based on dominance relation,  where they assign a fitness value (order/rank) for each individual in the population according to different properties based on the dominance relation, mainly: dominance rank, dominance count, and dominance depth \cite{2}.  The common point among these operators is that they can preserve Pareto solutions throughout the generational process. Furthermore, the current non-dominated solutions will iteratively approach the true Pareto front of the target problem. On the other hand, there are also criterion based selection operators that differ from the above mentioned ones, where it basically uses each of the objective functions periodically to select the solutions that pass to the next generation. In other words, preponderance is given to one of the objectives in each iteration of the algorithm, a classic example of criterion based techniques is VEGA  \cite{8}. In this work, we use both dominance based and criterion based selection techniques, where we extract and archive the non-dominated individuals, 	and then we assign fitness values to the rest of individuals according to their fitness in the $i^{th}$ criterion environment (only the $i^{th}$ objective function is considered, see Figure \ref{fig: eai}). This will keep the search process of each $MOEA_i$ focused on the $i^{th}$ criterion, while storing the non-dominated solutions met throughout the search process. Hence, the $MOEA_i$'s population will incrementally adapt to the $i^{th}$ criterion environment.   

For the selection operator of each $MOEA_i$, the order relation is defined as to give preponderance to the $i^{th}$ objective function, while retaining non-dominated solutions to be archived. The order relation is defined as follows. Let $PS_t$ be the set of  Pareto solutions obtained at iteration~$t$, 
$$ \forall x,y \in P_t, \; x\geq_i y \iff (x \succ y)\lor(Z^i(x)\geq Z^i(y)).\vspace{3mm}$$
Moreover, the number of dominated solution maintained for the next generation is determined by the population size $N$. More precisely, the parameter $N$ denotes the size of dominated individuals set, which are to be retained. Hence, a dominated solution is retained only if its rank according to the $i^{th}$ criterion is less than $N$, while the archive size is dynamic. The process of selecting individuals that pass to the next generation is given explicitly as follows: let $PS_{t}\subset P_t$, the set of current non-dominated individuals in $P_t$,  $PS_{t} = \{x\in P_t|\not\exists y\in P_t \, :\, y \succ x\}$. The population of the next generation $P_{t+1}$ is constructed as follows:
\begin{small}
\begin{eqnarray}
P_{t+1} &=&\left\{ x\in P_t \, | \, (x\in PS_{t})\lor \left(rank_i(x, P_t)\leq N + |PS_{t}|\right)\right\},\footnotemark\nonumber\\
 &=&\left\{ x\in P_t \, | \, (x\in PS_{t})\lor \left(rank_i(x, P_t \backslash PS_{t})\leq N\right) \right\}.\nonumber
\end{eqnarray}
\end{small}
\footnotetext{\ $rank_i(x,A)$ is the order of $x$ compared to elements of a set $A$, according to the $i^{th}$ objective function $Z^i$.}
In terms of complexity, the selection operator can be resumed in two major tasks: extracting non-dominated solutions running in $O\left(k|P_t|^2\right)$, and sorting dominated ones, which runs in the very worst case in $O\left(|P_t|\log(|P_t|)\right)$. Hence, the run-time of the selection procedure for each $MOEA_i$ is dominated by the process of extracting non-dominated solutions. Furthermore, the worst case complexity of the selection operator is $O\left(k|P_t|^2)\right)$.
\subsection{\textbf{Master entity (global selection and partitionning technique)}}\label{master}
The master entity preforms a periodical collection of the current archives from the criterion based MOEAs, in order to extract from them current global Pareto solutions; then partitions these potentially efficient solutions into $k$ subsets, where $k$ is the number of objective functions. The obtained subsets are used to update archives (current non-dominated solutions) of each search entity $MOEA_i$ (dedicated to the $i^{th}$ objective). The main goad of the partition procedure is the keep each worker focused on a given area, by adjusting periodically elite solutions of search entities according to the topology of the current Pareto front. In other words, it helps intensifying the search in each area assigned to parallel MOEAs. The partitioning procedure basically uses the statistical quantile of a given order $\alpha$ \cite{10}. Let $X$ a discrete random variable, a value $x_\alpha$ is called a quantile of order $\alpha$, $0 < \alpha < 1$, if:
$$ \mathbb{P}\left(X < x_\alpha\right) \leq \alpha \leq \mathbb{P}\left(X \leq x_\alpha\right).\vspace{3mm}$$ 
If we project this concept on the distribution of Pareto solutions in the objective space, the quantile of order $\alpha$ can be defined as the evaluation of the individual which divides the Pareto front into two parts such that: $\alpha$ proportion of the total number of individuals in the Pareto front are less than or equal to this quantile, and $(1-\alpha)$ proportion of the total number of individuals in the current Pareto front $PF_t$ (iteration $t$) are greater than this quantile. Furthermore, let $p_\alpha \in PS_t$ a non-dominated individual, $Z^i(p_\alpha)\in PF_t$ is said to be the Pareto front's quantile of order $\alpha$ according to the $i^{th}$ objective function, if:
$$\forall\,p \in PS_t, \quad \mathbb{P}\left( Z^i(p) < Z^i(p_\alpha) \right) \leq \alpha \leq \mathbb{P}\left( Z^i(p) \leq Z^i(p_\alpha) \right).\vspace{3mm}$$ 
Hence, for a given $\alpha$, we can construct for each objective function $Z^i, i\in \{1, \dots, k\}$, a partition $\{PS_t\backslash P^i, P^i\}$ using $p_\alpha$, the quantile of order $\alpha$, such that: 
$$P^i = \{ p\in PS_t| Z^i(p)\geq Z^i(p_{\alpha}) \}. $$
However, for a problem with $k$ objective functions, the set of non-dominated solutions is partitioned into $k$ subsets $\{P^i\}_{1\leq i \leq k}$ using the partitioning procedure (see Algorithm 2), where $PS_t = \cup_{1\leq i \leq k} P^i$, and each pair of the $\{P^i\}_{1\leq i \leq k}$ shares a subset of solution $F$. The size of $F$ is determined by the parameter $\alpha$, for example $|F|=(1-2\alpha)|PS_t|$ for bi-objective problems.  

\begin{figure}[H]
\label{alg:part}
\begin{flushleft}
\rule{\linewidth}{0.5pt}\vspace{1mm}\\
\textbf{Algorithm 2} Partitioning procedure\vspace{-2mm}\\  
\rule{\linewidth}{0.5pt}\vspace{0.1cm}
\textbf{Requires:} Set of non-dominated solutions $P$, and $\alpha \in\; ]0,0.5]$.\\
\textbf{Ensures:} $k$ subsets $\{P^i\}_{1\leq i \leq k}$ of a given set $P$.\vspace{-2mm}\\
\rule{\linewidth}{0.5pt}\vspace{0.1cm}
\hspace{5mm} \textbf{For} each objective $Z^i$ \textbf{do}\vspace{1mm}\\
\hspace{9mm} $P \leftarrow $ $sort(P, Z^i, \uparrow)$;\footnotemark
\vspace{1mm}\\
\hspace{9mm} $P^i \leftarrow  P \left( [\alpha|P|]:|P| \right)$;\footnotemark
\vspace{1mm}\\
\hspace{9mm} Send $P^i$ to worker $MOEA_i$;{\color {darkgray!90}\ // update $MOEA_i$ archive.}\vspace{1mm}\\
\textbf{End}\vspace{-1mm}\\
\rule{\linewidth}{0.5pt}\vspace{-3mm}\\
\end{flushleft}
\end{figure}

\footnotetext[3]{\ $sort(P, Z^i, \uparrow)$ denotes the procedure of sorting individuals in $P$ according to $Z^i$, in an increasing order.}
\footnotetext[4]{\ $P \left( a:b \right) = \{ p_j \in P | a\leq j \leq b \}$, and $\left[ a\right]$ is the integer part of $a$.}
  
Since the partitioning procedure is invoked periodically throughout the search process, and also needs to take all objective functions into consideration, it is imperative that the chosen partitioning procedure runs in a reasonable time.  Hence, we employed the concept of quantiles as to efficiently construct an adequate partition of the current non-dominated individuals recovered by the master. Whilst the partitioning is accomplished according to the distribution of non-dominated solutions in the decision space. In other words, all the objective functions are taken into consideration by the partitioning procedure.  The run-time of the partitioning procedure is $O\left(k |P| \log(|P|)\right)$. 
 The following algorithm resumes the functioning of the master entity. 

\begin{figure}[H]
\label{alg:part}
\begin{flushleft}
\rule{\linewidth}{0.5pt}\vspace{1mm}\\
\textbf{Algorithm 3} Master process\vspace{-2mm}\\  
\rule{\linewidth}{0.5pt}\vspace{0.1cm}
 \textbf{Repeat} \vspace{1mm}\\
\hspace{3mm}  \textbf{if} (Migration condition is satisfied) \textbf{then}
\vspace{1mm}\\
\hspace{7mm} Collect current archives $P^i$ from $MOEA_i, i=1,\dots, k$, $P=P\cup(\cup_i P^i)$;
\vspace{1mm}\\
\hspace{8mm}Remove dominated solutions from $P$;
\vspace{1mm}\\
\hspace{8mm}Preform partitioning procedure, $(P^1,\dots,P^k)\leftarrow Partition(P,\alpha)$;
\vspace{1mm}\\
\hspace{6mm} {\color {darkgray!90}// $P^1,\dots,P^k$ are sent to $MOEA_1, \dots, MOEA_k$, respectively.}\vspace{1mm}\\
\hspace{3mm}  \textbf{endif}\vspace{1mm}\\
 \textbf{until} (Termination condition met) \vspace{1mm}\\
\rule{\linewidth}{0.5pt}\vspace{-3mm}\\
\end{flushleft}
\end{figure}

\begin{figure}[h!]
\pgfplotsset{width=9.4cm,compat=newest, every tick label/.append style={font=\scriptsize}}
\begin{tikzpicture} [transform shape, scale = .9]
\begin{axis}[
        x=.165mm,
    legend cell align={left},  legend style={cells={align=left}},legend pos=south west,
]
\addplot[only marks, red!70,
    mark=square,mark=x]
    coordinates {
 (2845, 3001)  (2842, 3002)  (2841, 3005)  (2837, 3008)  (2834, 3013)  (2833, 3018)  (2831, 3021)  (2827, 3023)  (2825, 3029)  (2822, 3032)  (2819, 3037)  (2818, 3039)  (2812, 3042)  (2811, 3047)  (2807, 3049)  (2804, 3050)  (2803, 3052)  (2800, 3057)  (2800, 3057)  (2798, 3059)  (2797, 3065)  (2796, 3067)  (2793, 3068)  (2791, 3073)  (2791, 3073)  (2786, 3074)  (2785, 3075)  (2783, 3080)  (2782, 3084)  (2778, 3086)  (2778, 3086)  (2777, 3092)  (2773, 3095)  (2769, 3099)  (2766, 3102)  (2762, 3109)  (2756, 3110)  (2755, 3117)  (2751, 3119)  (2748, 3128)  (2742, 3129)  (2741, 3130)  (2740, 3133)  (2738, 3135)  (2737, 3137)  (2735, 3143)  (2729, 3144)  (2727, 3145)  (2723, 3149)  (2722, 3151)  (2717, 3155)  (2712, 3159)  (2712, 3159)  (2711, 3162)  (2707, 3168)  (2705, 3170)  (2704, 3176)  (2696, 3182)  (2690, 3191)  (2680, 3196)  (2674, 3200)  (2670, 3206)  (2662, 3207)  (2660, 3210)  (2655, 3213)  (2654, 3222)  (2647, 3223)  (2645, 3224)  (2643, 3226)  (2635, 3235)  (2630, 3236)  (2627, 3239)  (2623, 3242)  (2605, 3250)  (2602, 3253)  (2596, 3264)  (2591, 3266)  (2569, 3272)  (2564, 3277)  (2564, 3277)  (2556, 3278)  (2550, 3279)  (2549, 3284)  (2544, 3285)  (2538, 3286)  (2531, 3287)  (2529, 3290)  (2523, 3296)  (2519, 3299)  (2514, 3300)  (2510, 3301)  (2503, 3302)  (2498, 3303)  (2490, 3304)  (2484, 3312)  (2468, 3315)  (2442, 3316)  (2434, 3317)  (2431, 3321)  (2423, 3325)  (2408, 3327)  (2392, 3330)  (2375, 3334)  (2368, 3336)  (2328, 3339)  (2305, 3343)  (2277, 3344)  };
       \addlegendentry{\ {\scriptsize handled only by $MOEA_2$}}
\addplot[only marks, blue!70,
    mark=square,mark=x]
    coordinates {
 (2627, 3239)  (2630, 3236)  (2635, 3235)  (2643, 3226)  (2645, 3224)  (2647, 3223)  (2654, 3222)  (2655, 3213)  (2660, 3210)  (2662, 3207)  (2670, 3206)  (2674, 3200)  (2680, 3196)  (2690, 3191)  (2696, 3182)  (2704, 3176)  (2705, 3170)  (2707, 3168)  (2711, 3162)  (2712, 3159)  (2712, 3159)  (2717, 3155)  (2722, 3151)  (2723, 3149)  (2727, 3145)  (2729, 3144)  (2735, 3143)  (2737, 3137)  (2738, 3135)  (2740, 3133)  (2741, 3130)  (2742, 3129)  (2748, 3128)  (2751, 3119)  (2755, 3117)  (2756, 3110)  (2762, 3109)  (2766, 3102)  (2769, 3099)  (2773, 3095)  (2777, 3092)  (2778, 3086)  (2778, 3086)  (2782, 3084)  (2783, 3080)  (2785, 3075)  (2786, 3074)  (2791, 3073)  (2791, 3073)  (2793, 3068)  (2796, 3067)  (2797, 3065)  (2798, 3059)  (2800, 3057)  (2800, 3057)  (2803, 3052)  (2804, 3050)  (2807, 3049)  (2811, 3047)  (2812, 3042)  (2818, 3039)  (2819, 3037)  (2822, 3032)  (2825, 3029)  (2827, 3023)  (2831, 3021)  (2833, 3018)  (2834, 3013)  (2837, 3008)  (2841, 3005)  (2842, 3002)  (2845, 3001)  (2846, 2992)  (2846, 2992)  (2850, 2989)  (2853, 2985)  (2857, 2975)  (2864, 2970)  (2872, 2954)  (2875, 2938)  (2877, 2937)  (2878, 2933)  (2881, 2929)  (2885, 2920)  (2892, 2916)  (2894, 2904)  (2900, 2898)  (2902, 2882)  (2904, 2877)  (2904, 2877)  (2906, 2874)  (2907, 2871)  (2909, 2864)  (2913, 2861)  (2915, 2849)  (2917, 2837)  (2919, 2836)  (2923, 2831)  (2925, 2827)  (2926, 2822)  (2934, 2815)  (2936, 2788)  (2938, 2783)  (2940, 2780)  (2945, 2728)  (2946, 2698)  (2951, 2651)  (2627, 3239)  (2630, 3236)  (2635, 3235)  (2643, 3226)  (2645, 3224)  (2647, 3223)  (2654, 3222)  (2655, 3213)  (2660, 3210)  (2662, 3207)  (2670, 3206)  (2674, 3200)  (2680, 3196)  (2690, 3191)  (2696, 3182)  (2704, 3176)  (2705, 3170)  (2707, 3168)  (2711, 3162)  (2712, 3159)  (2712, 3159)  (2717, 3155)  (2722, 3151)  (2723, 3149)  (2727, 3145)  (2729, 3144)  (2735, 3143)  (2737, 3137)  (2738, 3135)  (2740, 3133)  (2741, 3130)  (2742, 3129)  (2748, 3128)  (2751, 3119)  (2755, 3117)  (2756, 3110)  (2762, 3109)  (2766, 3102)  (2769, 3099)  (2773, 3095)  (2777, 3092)  (2778, 3086)  (2778, 3086)  (2782, 3084)  (2783, 3080)  (2785, 3075)  (2786, 3074)  (2791, 3073)  (2791, 3073)  (2793, 3068)  (2796, 3067)  (2797, 3065)  (2798, 3059)  (2800, 3057)  (2800, 3057)  (2803, 3052)  (2804, 3050)  (2807, 3049)  (2811, 3047)  (2812, 3042)  (2818, 3039)  (2819, 3037)  (2822, 3032)  (2825, 3029)  (2827, 3023)  (2831, 3021)  (2833, 3018)  (2834, 3013)  (2837, 3008)  (2841, 3005)  (2842, 3002)  (2845, 3001)  (2846, 2992)  (2846, 2992)  (2850, 2989)  (2853, 2985)  (2857, 2975)  (2864, 2970)  (2872, 2954)  (2875, 2938)  (2877, 2937)  (2878, 2933)  (2881, 2929)  (2885, 2920)  (2892, 2916)  (2894, 2904)  (2900, 2898)  (2902, 2882)  (2904, 2877)  (2904, 2877)  (2906, 2874)  (2907, 2871)  (2909, 2864)  (2913, 2861)  (2915, 2849)  (2917, 2837)  (2919, 2836)  (2923, 2831)  (2925, 2827)  (2926, 2822)  (2934, 2815)  (2936, 2788)  (2938, 2783)  (2940, 2780)  (2945, 2728)  (2946, 2698)  (2951, 2651)  (2627, 3239)  (2630, 3236)  (2635, 3235)  (2643, 3226)  (2645, 3224)  (2647, 3223)  (2654, 3222)  (2655, 3213)  (2660, 3210)  (2662, 3207)  (2670, 3206)  (2674, 3200)  (2680, 3196)  (2690, 3191)  (2696, 3182)  (2704, 3176)  (2705, 3170)  (2707, 3168)  (2711, 3162)  (2712, 3159)  (2712, 3159)  (2717, 3155)  (2722, 3151)  (2723, 3149)  (2727, 3145)  (2729, 3144)  (2735, 3143)  (2737, 3137)  (2738, 3135)  (2740, 3133)  (2741, 3130)  (2742, 3129)  (2748, 3128)  (2751, 3119)  (2755, 3117)  (2756, 3110)  (2762, 3109)  (2766, 3102)  (2769, 3099)  (2773, 3095)  (2777, 3092)  (2778, 3086)  (2778, 3086)  (2782, 3084)  (2783, 3080)  (2785, 3075)  (2786, 3074)  (2791, 3073)  (2791, 3073)  (2793, 3068)  (2796, 3067)  (2797, 3065)  (2798, 3059)  (2800, 3057)  (2800, 3057)  (2803, 3052)  (2804, 3050)  (2807, 3049)  (2811, 3047)  (2812, 3042)  (2818, 3039)  (2819, 3037)  (2822, 3032)  (2825, 3029)  (2827, 3023)  (2831, 3021)  (2833, 3018)  (2834, 3013)  (2837, 3008)  (2841, 3005)  (2842, 3002)  (2845, 3001)  (2846, 2992)  (2846, 2992)  (2850, 2989)  (2853, 2985)  (2857, 2975)  (2864, 2970)  (2872, 2954)  (2875, 2938)  (2877, 2937)  (2878, 2933)  (2881, 2929)  (2885, 2920)  (2892, 2916)  (2894, 2904)  (2900, 2898)  (2902, 2882)  (2904, 2877)  (2904, 2877)  (2906, 2874)  (2907, 2871)  (2909, 2864)  (2913, 2861)  (2915, 2849)  (2917, 2837)  (2919, 2836)  (2923, 2831)  (2925, 2827)  (2926, 2822)  (2934, 2815)  (2936, 2788)  (2938, 2783)  (2940, 2780)  (2945, 2728)  (2946, 2698)  (2951, 2651)  (2627, 3239)  (2630, 3236)  (2635, 3235)  (2643, 3226)  (2645, 3224)  (2647, 3223)  (2654, 3222)  (2655, 3213)  (2660, 3210)  (2662, 3207)  (2670, 3206)  (2674, 3200)  (2680, 3196)  (2690, 3191)  (2696, 3182)  (2704, 3176)  (2705, 3170)  (2707, 3168)  (2711, 3162)  (2712, 3159)  (2712, 3159)  (2717, 3155)  (2722, 3151)  (2723, 3149)  (2727, 3145)  (2729, 3144)  (2735, 3143)  (2737, 3137)  (2738, 3135)  (2740, 3133)  (2741, 3130)  (2742, 3129)  (2748, 3128)  (2751, 3119)  (2755, 3117)  (2756, 3110)  (2762, 3109)  (2766, 3102)  (2769, 3099)  (2773, 3095)  (2777, 3092)  (2778, 3086)  (2778, 3086)  (2782, 3084)  (2783, 3080)  (2785, 3075)  (2786, 3074)  (2791, 3073)  (2791, 3073)  (2793, 3068)  (2796, 3067)  (2797, 3065)  (2798, 3059)  (2800, 3057)  (2800, 3057)  (2803, 3052)  (2804, 3050)  (2807, 3049)  (2811, 3047)  (2812, 3042)  (2818, 3039)  (2819, 3037)  (2822, 3032)  (2825, 3029)  (2827, 3023)  (2831, 3021)  (2833, 3018)  (2834, 3013)  (2837, 3008)  (2841, 3005)  (2842, 3002)  (2845, 3001)  (2846, 2992)  (2846, 2992)  (2850, 2989)  (2853, 2985)  (2857, 2975)  (2864, 2970)  (2872, 2954)  (2875, 2938)  (2877, 2937)  (2878, 2933)  (2881, 2929)  (2885, 2920)  (2892, 2916)  (2894, 2904)  (2900, 2898)  (2902, 2882)  (2904, 2877)  (2904, 2877)  (2906, 2874)  (2907, 2871)  (2909, 2864)  (2913, 2861)  (2915, 2849)  (2917, 2837)  (2919, 2836)  (2923, 2831)  (2925, 2827)  (2926, 2822)  (2934, 2815)  (2936, 2788)  (2938, 2783)  (2940, 2780)  (2945, 2728)  (2946, 2698)  (2951, 2651)  (2627, 3239)  (2630, 3236)  (2635, 3235)  (2643, 3226)  (2645, 3224)  (2647, 3223)  (2654, 3222)  (2655, 3213)  (2660, 3210)  (2662, 3207)  (2670, 3206)  (2674, 3200)  (2680, 3196)  (2690, 3191)  (2696, 3182)  (2704, 3176)  (2705, 3170)  (2707, 3168)  (2711, 3162)  (2712, 3159)  (2712, 3159)  (2717, 3155)  (2722, 3151)  (2723, 3149)  (2727, 3145)  (2729, 3144)  (2735, 3143)  (2737, 3137)  (2738, 3135)  (2740, 3133)  (2741, 3130)  (2742, 3129)  (2748, 3128)  (2751, 3119)  (2755, 3117)  (2756, 3110)  (2762, 3109)  (2766, 3102)  (2769, 3099)  (2773, 3095)  (2777, 3092)  (2778, 3086)  (2778, 3086)  (2782, 3084)  (2783, 3080)  (2785, 3075)  (2786, 3074)  (2791, 3073)  (2791, 3073)  (2793, 3068)  (2796, 3067)  (2797, 3065)  (2798, 3059)  (2800, 3057)  (2800, 3057)  (2803, 3052)  (2804, 3050)  (2807, 3049)  (2811, 3047)  (2812, 3042)  (2818, 3039)  (2819, 3037)  (2822, 3032)  (2825, 3029)  (2827, 3023)  (2831, 3021)  (2833, 3018)  (2834, 3013)  (2837, 3008)  (2841, 3005)  (2842, 3002)  (2845, 3001)  (2846, 2992)  (2846, 2992)  (2850, 2989)  (2853, 2985)  (2857, 2975)  (2864, 2970)  (2872, 2954)  (2875, 2938)  (2877, 2937)  (2878, 2933)  (2881, 2929)  (2885, 2920)  (2892, 2916)  (2894, 2904)  (2900, 2898)  (2902, 2882)  (2904, 2877)  (2904, 2877)  (2906, 2874)  (2907, 2871)  (2909, 2864)  (2913, 2861)  (2915, 2849)  (2917, 2837)  (2919, 2836)  (2923, 2831)  (2925, 2827)  (2926, 2822)  (2934, 2815)  (2936, 2788)  (2938, 2783)  (2940, 2780)  (2945, 2728)  (2946, 2698)  (2951, 2651)  (2627, 3239)  (2630, 3236)  (2635, 3235)  (2643, 3226)  (2645, 3224)  (2647, 3223)  (2654, 3222)  (2655, 3213)  (2660, 3210)  (2662, 3207)  (2670, 3206)  (2674, 3200)  (2680, 3196)  (2690, 3191)  (2696, 3182)  (2704, 3176)  (2705, 3170)  (2707, 3168)  (2711, 3162)  (2712, 3159)  (2712, 3159)  (2717, 3155)  (2722, 3151)  (2723, 3149)  (2727, 3145)  (2729, 3144)  (2735, 3143)  (2737, 3137)  (2738, 3135)  (2740, 3133)  (2741, 3130)  (2742, 3129)  (2748, 3128)  (2751, 3119)  (2755, 3117)  (2756, 3110)  (2762, 3109)  (2766, 3102)  (2769, 3099)  (2773, 3095)  (2777, 3092)  (2778, 3086)  (2778, 3086)  (2782, 3084)  (2783, 3080)  (2785, 3075)  (2786, 3074)  (2791, 3073)  (2791, 3073)  (2793, 3068)  (2796, 3067)  (2797, 3065)  (2798, 3059)  (2800, 3057)  (2800, 3057)  (2803, 3052)  (2804, 3050)  (2807, 3049)  (2811, 3047)  (2812, 3042)  (2818, 3039)  (2819, 3037)  (2822, 3032)  (2825, 3029)  (2827, 3023)  (2831, 3021)  (2833, 3018)  (2834, 3013)  (2837, 3008)  (2841, 3005)  (2842, 3002)  (2845, 3001)  (2846, 2992)  (2846, 2992)  (2850, 2989)  (2853, 2985)  (2857, 2975)  (2864, 2970)  (2872, 2954)  (2875, 2938)  (2877, 2937)  (2878, 2933)  (2881, 2929)  (2885, 2920)  (2892, 2916)  (2894, 2904)  (2900, 2898)  (2902, 2882)  (2904, 2877)  (2904, 2877)  (2906, 2874)  (2907, 2871)  (2909, 2864)  (2913, 2861)  (2915, 2849)  (2917, 2837)  (2919, 2836)  (2923, 2831)  (2925, 2827)  (2926, 2822)  (2934, 2815)  (2936, 2788)  (2938, 2783)  (2940, 2780)  (2945, 2728)  (2946, 2698)  (2951, 2651)  (2627, 3239)  (2630, 3236)  (2635, 3235)  (2643, 3226)  (2645, 3224)  (2647, 3223)  (2654, 3222)  (2655, 3213)  (2660, 3210)  (2662, 3207)  (2670, 3206)  (2674, 3200)  (2680, 3196)  (2690, 3191)  (2696, 3182)  (2704, 3176)  (2705, 3170)  (2707, 3168)  (2711, 3162)  (2712, 3159)  (2712, 3159)  (2717, 3155)  (2722, 3151)  (2723, 3149)  (2727, 3145)  (2729, 3144)  (2735, 3143)  (2737, 3137)  (2738, 3135)  (2740, 3133)  (2741, 3130)  (2742, 3129)  (2748, 3128)  (2751, 3119)  (2755, 3117)  (2756, 3110)  (2762, 3109)  (2766, 3102)  (2769, 3099)  (2773, 3095)  (2777, 3092)  (2778, 3086)  (2778, 3086)  (2782, 3084)  (2783, 3080)  (2785, 3075)  (2786, 3074)  (2791, 3073)  (2791, 3073)  (2793, 3068)  (2796, 3067)  (2797, 3065)  (2798, 3059)  (2800, 3057)  (2800, 3057)  (2803, 3052)  (2804, 3050)  (2807, 3049)  (2811, 3047)  (2812, 3042)  (2818, 3039)  (2819, 3037)  (2822, 3032)  (2825, 3029)  (2827, 3023)  (2831, 3021)  (2833, 3018)  (2834, 3013)  (2837, 3008)  (2841, 3005)  (2842, 3002)  (2845, 3001)  (2846, 2992)  (2846, 2992)  (2850, 2989)  (2853, 2985)  (2857, 2975)  (2864, 2970)  (2872, 2954)  (2875, 2938)  (2877, 2937)  (2878, 2933)  (2881, 2929)  (2885, 2920)  (2892, 2916)  (2894, 2904)  (2900, 2898)  (2902, 2882)  (2904, 2877)  (2904, 2877)  (2906, 2874)  (2907, 2871)  (2909, 2864)  (2913, 2861)  (2915, 2849)  (2917, 2837)  (2919, 2836)  (2923, 2831)  (2925, 2827)  (2926, 2822)  (2934, 2815)  (2936, 2788)  (2938, 2783)  (2940, 2780)  (2945, 2728)  (2946, 2698)  (2951, 2651)  (2627, 3239)  (2630, 3236)  (2635, 3235)  (2643, 3226)  (2645, 3224)  (2647, 3223)  (2654, 3222)  (2655, 3213)  (2660, 3210)  (2662, 3207)  (2670, 3206)  (2674, 3200)  (2680, 3196)  (2690, 3191)  (2696, 3182)  (2704, 3176)  (2705, 3170)  (2707, 3168)  (2711, 3162)  (2712, 3159)  (2712, 3159)  (2717, 3155)  (2722, 3151)  (2723, 3149)  (2727, 3145)  (2729, 3144)  (2735, 3143)  (2737, 3137)  (2738, 3135)  (2740, 3133)  (2741, 3130)  (2742, 3129)  (2748, 3128)  (2751, 3119)  (2755, 3117)  (2756, 3110)  (2762, 3109)  (2766, 3102)  (2769, 3099)  (2773, 3095)  (2777, 3092)  (2778, 3086)  (2778, 3086)  (2782, 3084)  (2783, 3080)  (2785, 3075)  (2786, 3074)  (2791, 3073)  (2791, 3073)  (2793, 3068)  (2796, 3067)  (2797, 3065)  (2798, 3059)  (2800, 3057)  (2800, 3057)  (2803, 3052)  (2804, 3050)  (2807, 3049)  (2811, 3047)  (2812, 3042)  (2818, 3039)  (2819, 3037)  (2822, 3032)  (2825, 3029)  (2827, 3023)  (2831, 3021)  (2833, 3018)  (2834, 3013)  (2837, 3008)  (2841, 3005)  (2842, 3002)  (2845, 3001)  (2846, 2992)  (2846, 2992)  (2850, 2989)  (2853, 2985)  (2857, 2975)  (2864, 2970)  (2872, 2954)  (2875, 2938)  (2877, 2937)  (2878, 2933)  (2881, 2929)  (2885, 2920)  (2892, 2916)  (2894, 2904)  (2900, 2898)  (2902, 2882)  (2904, 2877)  (2904, 2877)  (2906, 2874)  (2907, 2871)  (2909, 2864)  (2913, 2861)  (2915, 2849)  (2917, 2837)  (2919, 2836)  (2923, 2831)  (2925, 2827)  (2926, 2822)  (2934, 2815)  (2936, 2788)  (2938, 2783)  (2940, 2780)  (2945, 2728)  (2946, 2698)  (2951, 2651)  (2627, 3239)  (2630, 3236)  (2635, 3235)  (2643, 3226)  (2645, 3224)  (2647, 3223)  (2654, 3222)  (2655, 3213)  (2660, 3210)  (2662, 3207)  (2670, 3206)  (2674, 3200)  (2680, 3196)  (2690, 3191)  (2696, 3182)  (2704, 3176)  (2705, 3170)  (2707, 3168)  (2711, 3162)  (2712, 3159)  (2712, 3159)  (2717, 3155)  (2722, 3151)  (2723, 3149)  (2727, 3145)  (2729, 3144)  (2735, 3143)  (2737, 3137)  (2738, 3135)  (2740, 3133)  (2741, 3130)  (2742, 3129)  (2748, 3128)  (2751, 3119)  (2755, 3117)  (2756, 3110)  (2762, 3109)  (2766, 3102)  (2769, 3099)  (2773, 3095)  (2777, 3092)  (2778, 3086)  (2778, 3086)  (2782, 3084)  (2783, 3080)  (2785, 3075)  (2786, 3074)  (2791, 3073)  (2791, 3073)  (2793, 3068)  (2796, 3067)  (2797, 3065)  (2798, 3059)  (2800, 3057)  (2800, 3057)  (2803, 3052)  (2804, 3050)  (2807, 3049)  (2811, 3047)  (2812, 3042)  (2818, 3039)  (2819, 3037)  (2822, 3032)  (2825, 3029)  (2827, 3023)  (2831, 3021)  (2833, 3018)  (2834, 3013)  (2837, 3008)  (2841, 3005)  (2842, 3002)  (2845, 3001)  (2846, 2992)  (2846, 2992)  (2850, 2989)  (2853, 2985)  (2857, 2975)  (2864, 2970)  (2872, 2954)  (2875, 2938)  (2877, 2937)  (2878, 2933)  (2881, 2929)  (2885, 2920)  (2892, 2916)  (2894, 2904)  (2900, 2898)  (2902, 2882)  (2904, 2877)  (2904, 2877)  (2906, 2874)  (2907, 2871)  (2909, 2864)  (2913, 2861)  (2915, 2849)  (2917, 2837)  (2919, 2836)  (2923, 2831)  (2925, 2827)  (2926, 2822)  (2934, 2815)  (2936, 2788)  (2938, 2783)  (2940, 2780)  (2945, 2728)  (2946, 2698)  (2951, 2651)  (2627, 3239)  (2630, 3236)  (2635, 3235)  (2643, 3226)  (2645, 3224)  (2647, 3223)  (2654, 3222)  (2655, 3213)  (2660, 3210)  (2662, 3207)  (2670, 3206)  (2674, 3200)  (2680, 3196)  (2690, 3191)  (2696, 3182)  (2704, 3176)  (2705, 3170)  (2707, 3168)  (2711, 3162)  (2712, 3159)  (2712, 3159)  (2717, 3155)  (2722, 3151)  (2723, 3149)  (2727, 3145)  (2729, 3144)  (2735, 3143)  (2737, 3137)  (2738, 3135)  (2740, 3133)  (2741, 3130)  (2742, 3129)  (2748, 3128)  (2751, 3119)  (2755, 3117)  (2756, 3110)  (2762, 3109)  (2766, 3102)  (2769, 3099)  (2773, 3095)  (2777, 3092)  (2778, 3086)  (2778, 3086)  (2782, 3084)  (2783, 3080)  (2785, 3075)  (2786, 3074)  (2791, 3073)  (2791, 3073)  (2793, 3068)  (2796, 3067)  (2797, 3065)  (2798, 3059)  (2800, 3057)  (2800, 3057)  (2803, 3052)  (2804, 3050)  (2807, 3049)  (2811, 3047)  (2812, 3042)  (2818, 3039)  (2819, 3037)  (2822, 3032)  (2825, 3029)  (2827, 3023)  (2831, 3021)  (2833, 3018)  (2834, 3013)  (2837, 3008)  (2841, 3005)  (2842, 3002)  (2845, 3001)  (2846, 2992)  (2846, 2992)  (2850, 2989)  (2853, 2985)  (2857, 2975)  (2864, 2970)  (2872, 2954)  (2875, 2938)  (2877, 2937)  (2878, 2933)  (2881, 2929)  (2885, 2920)  (2892, 2916)  (2894, 2904)  (2900, 2898)  (2902, 2882)  (2904, 2877)  (2904, 2877)  (2906, 2874)  (2907, 2871)  (2909, 2864)  (2913, 2861)  (2915, 2849)  (2917, 2837)  (2919, 2836)  (2923, 2831)  (2925, 2827)  (2926, 2822)  (2934, 2815)  (2936, 2788)  (2938, 2783)  (2940, 2780)  (2945, 2728)  (2946, 2698)  (2951, 2651)  (2627, 3239)  (2630, 3236)  (2635, 3235)  (2643, 3226)  (2645, 3224)  (2647, 3223)  (2654, 3222)  (2655, 3213)  (2660, 3210)  (2662, 3207)  (2670, 3206)  (2674, 3200)  (2680, 3196)  (2690, 3191)  (2696, 3182)  (2704, 3176)  (2705, 3170)  (2707, 3168)  (2711, 3162)  (2712, 3159)  (2712, 3159)  (2717, 3155)  (2722, 3151)  (2723, 3149)  (2727, 3145)  (2729, 3144)  (2735, 3143)  (2737, 3137)  (2738, 3135)  (2740, 3133)  (2741, 3130)  (2742, 3129)  (2748, 3128)  (2751, 3119)  (2755, 3117)  (2756, 3110)  (2762, 3109)  (2766, 3102)  (2769, 3099)  (2773, 3095)  (2777, 3092)  (2778, 3086)  (2778, 3086)  (2782, 3084)  (2783, 3080)  (2785, 3075)  (2786, 3074)  (2791, 3073)  (2791, 3073)  (2793, 3068)  (2796, 3067)  (2797, 3065)  (2798, 3059)  (2800, 3057)  (2800, 3057)  (2803, 3052)  (2804, 3050)  (2807, 3049)  (2811, 3047)  (2812, 3042)  (2818, 3039)  (2819, 3037)  (2822, 3032)  (2825, 3029)  (2827, 3023)  (2831, 3021)  (2833, 3018)  (2834, 3013)  (2837, 3008)  (2841, 3005)  (2842, 3002)  (2845, 3001)  (2846, 2992)  (2846, 2992)  (2850, 2989)  (2853, 2985)  (2857, 2975)  (2864, 2970)  (2872, 2954)  (2875, 2938)  (2877, 2937)  (2878, 2933)  (2881, 2929)  (2885, 2920)  (2892, 2916)  (2894, 2904)  (2900, 2898)  (2902, 2882)  (2904, 2877)  (2904, 2877)  (2906, 2874)  (2907, 2871)  (2909, 2864)  (2913, 2861)  (2915, 2849)  (2917, 2837)  (2919, 2836)  (2923, 2831)  (2925, 2827)  (2926, 2822)  (2934, 2815)  (2936, 2788)  (2938, 2783)  (2940, 2780)  (2945, 2728)  (2946, 2698)  (2951, 2651)  (2627, 3239)  (2630, 3236)  (2635, 3235)  (2643, 3226)  (2645, 3224)  (2647, 3223)  (2654, 3222)  (2655, 3213)  (2660, 3210)  (2662, 3207)  (2670, 3206)  (2674, 3200)  (2680, 3196)  (2690, 3191)  (2696, 3182)  (2704, 3176)  (2705, 3170)  (2707, 3168)  (2711, 3162)  (2712, 3159)  (2712, 3159)  (2717, 3155)  (2722, 3151)  (2723, 3149)  (2727, 3145)  (2729, 3144)  (2735, 3143)  (2737, 3137)  (2738, 3135)  (2740, 3133)  (2741, 3130)  (2742, 3129)  (2748, 3128)  (2751, 3119)  (2755, 3117)  (2756, 3110)  (2762, 3109)  (2766, 3102)  (2769, 3099)  (2773, 3095)  (2777, 3092)  (2778, 3086)  (2778, 3086)  (2782, 3084)  (2783, 3080)  (2785, 3075)  (2786, 3074)  (2791, 3073)  (2791, 3073)  (2793, 3068)  (2796, 3067)  (2797, 3065)  (2798, 3059)  (2800, 3057)  (2800, 3057)  (2803, 3052)  (2804, 3050)  (2807, 3049)  (2811, 3047)  (2812, 3042)  (2818, 3039)  (2819, 3037)  (2822, 3032)  (2825, 3029)  (2827, 3023)  (2831, 3021)  (2833, 3018)  (2834, 3013)  (2837, 3008)  (2841, 3005)  (2842, 3002)  (2845, 3001)  (2846, 2992)  (2846, 2992)  (2850, 2989)  (2853, 2985)  (2857, 2975)  (2864, 2970)  (2872, 2954)  (2875, 2938)  (2877, 2937)  (2878, 2933)  (2881, 2929)  (2885, 2920)  (2892, 2916)  (2894, 2904)  (2900, 2898)  (2902, 2882)  (2904, 2877)  (2904, 2877)  (2906, 2874)  (2907, 2871)  (2909, 2864)  (2913, 2861)  (2915, 2849)  (2917, 2837)  (2919, 2836)  (2923, 2831)  (2925, 2827)  (2926, 2822)  (2934, 2815)  (2936, 2788)  (2938, 2783)  (2940, 2780)  (2945, 2728)  (2946, 2698)  (2951, 2651)  (2627, 3239)  (2630, 3236)  (2635, 3235)  (2643, 3226)  (2645, 3224)  (2647, 3223)  (2654, 3222)  (2655, 3213)  (2660, 3210)  (2662, 3207)  (2670, 3206)  (2674, 3200)  (2680, 3196)  (2690, 3191)  (2696, 3182)  (2704, 3176)  (2705, 3170)  (2707, 3168)  (2711, 3162)  (2712, 3159)  (2712, 3159)  (2717, 3155)  (2722, 3151)  (2723, 3149)  (2727, 3145)  (2729, 3144)  (2735, 3143)  (2737, 3137)  (2738, 3135)  (2740, 3133)  (2741, 3130)  (2742, 3129)  (2748, 3128)  (2751, 3119)  (2755, 3117)  (2756, 3110)  (2762, 3109)  (2766, 3102)  (2769, 3099)  (2773, 3095)  (2777, 3092)  (2778, 3086)  (2778, 3086)  (2782, 3084)  (2783, 3080)  (2785, 3075)  (2786, 3074)  (2791, 3073)  (2791, 3073)  (2793, 3068)  (2796, 3067)  (2797, 3065)  (2798, 3059)  (2800, 3057)  (2800, 3057)  (2803, 3052)  (2804, 3050)  (2807, 3049)  (2811, 3047)  (2812, 3042)  (2818, 3039)  (2819, 3037)  (2822, 3032)  (2825, 3029)  (2827, 3023)  (2831, 3021)  (2833, 3018)  (2834, 3013)  (2837, 3008)  (2841, 3005)  (2842, 3002)  (2845, 3001)  (2846, 2992)  (2846, 2992)  (2850, 2989)  (2853, 2985)  (2857, 2975)  (2864, 2970)  (2872, 2954)  (2875, 2938)  (2877, 2937)  (2878, 2933)  (2881, 2929)  (2885, 2920)  (2892, 2916)  (2894, 2904)  (2900, 2898)  (2902, 2882)  (2904, 2877)  (2904, 2877)  (2906, 2874)  (2907, 2871)  (2909, 2864)  (2913, 2861)  (2915, 2849)  (2917, 2837)  (2919, 2836)  (2923, 2831)  (2925, 2827)  (2926, 2822)  (2934, 2815)  (2936, 2788)  (2938, 2783)  (2940, 2780)  (2945, 2728)  (2946, 2698)  (2951, 2651)  };
      \addlegendentry{\ {\scriptsize handled only by $MOEA_1$}}
\addplot[only marks, green!70,
    mark=square,mark=x]
   coordinates {
 (2627, 3239)  (2630, 3236)  (2635, 3235)  (2643, 3226)  (2645, 3224)  (2647, 3223)  (2654, 3222)  (2655, 3213)  (2660, 3210)  (2662, 3207)  (2670, 3206)  (2674, 3200)  (2680, 3196)  (2690, 3191)  (2696, 3182)  (2704, 3176)  (2705, 3170)  (2707, 3168)  (2711, 3162)  (2712, 3159)  (2712, 3159)  (2717, 3155)  (2722, 3151)  (2723, 3149)  (2727, 3145)  (2729, 3144)  (2735, 3143)  (2737, 3137)  (2738, 3135)  (2740, 3133)  (2741, 3130)  (2742, 3129)  (2748, 3128)  (2751, 3119)  (2755, 3117)  (2756, 3110)  (2762, 3109)  (2766, 3102)  (2769, 3099)  (2773, 3095)  (2777, 3092)  (2778, 3086)  (2778, 3086)  (2782, 3084)  (2783, 3080)  (2785, 3075)  (2786, 3074)  (2791, 3073)  (2791, 3073)  (2793, 3068)  (2796, 3067)  (2797, 3065)  (2798, 3059)  (2800, 3057)  (2800, 3057)  (2803, 3052)  (2804, 3050)  (2807, 3049)  (2811, 3047)  (2812, 3042)  (2818, 3039)  (2819, 3037)  (2822, 3032)  (2825, 3029)  (2827, 3023)  (2831, 3021)  (2833, 3018)  (2834, 3013)  (2837, 3008)  (2841, 3005)  (2842, 3002)  (2845, 3001)  (2627, 3239)  (2630, 3236)  (2635, 3235)  (2643, 3226)  (2645, 3224)  (2647, 3223)  (2654, 3222)  (2655, 3213)  (2660, 3210)  (2662, 3207)  (2670, 3206)  (2674, 3200)  (2680, 3196)  (2690, 3191)  (2696, 3182)  (2704, 3176)  (2705, 3170)  (2707, 3168)  (2711, 3162)  (2712, 3159)  (2712, 3159)  (2717, 3155)  (2722, 3151)  (2723, 3149)  (2727, 3145)  (2729, 3144)  (2735, 3143)  (2737, 3137)  (2738, 3135)  (2740, 3133)  (2741, 3130)  (2742, 3129)  (2748, 3128)  (2751, 3119)  (2755, 3117)  (2756, 3110)  (2762, 3109)  (2766, 3102)  (2769, 3099)  (2773, 3095)  (2777, 3092)  (2778, 3086)  (2778, 3086)  (2782, 3084)  (2783, 3080)  (2785, 3075)  (2786, 3074)  (2791, 3073)  (2791, 3073)  (2793, 3068)  (2796, 3067)  (2797, 3065)  (2798, 3059)  (2800, 3057)  (2800, 3057)  (2803, 3052)  (2804, 3050)  (2807, 3049)  (2811, 3047)  (2812, 3042)  (2818, 3039)  (2819, 3037)  (2822, 3032)  (2825, 3029)  (2827, 3023)  (2831, 3021)  (2833, 3018)  (2834, 3013)  (2837, 3008)  (2841, 3005)  (2842, 3002)  (2845, 3001)  (2627, 3239)  (2630, 3236)  (2635, 3235)  (2643, 3226)  (2645, 3224)  (2647, 3223)  (2654, 3222)  (2655, 3213)  (2660, 3210)  (2662, 3207)  (2670, 3206)  (2674, 3200)  (2680, 3196)  (2690, 3191)  (2696, 3182)  (2704, 3176)  (2705, 3170)  (2707, 3168)  (2711, 3162)  (2712, 3159)  (2712, 3159)  (2717, 3155)  (2722, 3151)  (2723, 3149)  (2727, 3145)  (2729, 3144)  (2735, 3143)  (2737, 3137)  (2738, 3135)  (2740, 3133)  (2741, 3130)  (2742, 3129)  (2748, 3128)  (2751, 3119)  (2755, 3117)  (2756, 3110)  (2762, 3109)  (2766, 3102)  (2769, 3099)  (2773, 3095)  (2777, 3092)  (2778, 3086)  (2778, 3086)  (2782, 3084)  (2783, 3080)  (2785, 3075)  (2786, 3074)  (2791, 3073)  (2791, 3073)  (2793, 3068)  (2796, 3067)  (2797, 3065)  (2798, 3059)  (2800, 3057)  (2800, 3057)  (2803, 3052)  (2804, 3050)  (2807, 3049)  (2811, 3047)  (2812, 3042)  (2818, 3039)  (2819, 3037)  (2822, 3032)  (2825, 3029)  (2827, 3023)  (2831, 3021)  (2833, 3018)  (2834, 3013)  (2837, 3008)  (2841, 3005)  (2842, 3002)  (2845, 3001)  (2627, 3239)  (2630, 3236)  (2635, 3235)  (2643, 3226)  (2645, 3224)  (2647, 3223)  (2654, 3222)  (2655, 3213)  (2660, 3210)  (2662, 3207)  (2670, 3206)  (2674, 3200)  (2680, 3196)  (2690, 3191)  (2696, 3182)  (2704, 3176)  (2705, 3170)  (2707, 3168)  (2711, 3162)  (2712, 3159)  (2712, 3159)  (2717, 3155)  (2722, 3151)  (2723, 3149)  (2727, 3145)  (2729, 3144)  (2735, 3143)  (2737, 3137)  (2738, 3135)  (2740, 3133)  (2741, 3130)  (2742, 3129)  (2748, 3128)  (2751, 3119)  (2755, 3117)  (2756, 3110)  (2762, 3109)  (2766, 3102)  (2769, 3099)  (2773, 3095)  (2777, 3092)  (2778, 3086)  (2778, 3086)  (2782, 3084)  (2783, 3080)  (2785, 3075)  (2786, 3074)  (2791, 3073)  (2791, 3073)  (2793, 3068)  (2796, 3067)  (2797, 3065)  (2798, 3059)  (2800, 3057)  (2800, 3057)  (2803, 3052)  (2804, 3050)  (2807, 3049)  (2811, 3047)  (2812, 3042)  (2818, 3039)  (2819, 3037)  (2822, 3032)  (2825, 3029)  (2827, 3023)  (2831, 3021)  (2833, 3018)  (2834, 3013)  (2837, 3008)  (2841, 3005)  (2842, 3002)  (2845, 3001)  (2627, 3239)  (2630, 3236)  (2635, 3235)  (2643, 3226)  (2645, 3224)  (2647, 3223)  (2654, 3222)  (2655, 3213)  (2660, 3210)  (2662, 3207)  (2670, 3206)  (2674, 3200)  (2680, 3196)  (2690, 3191)  (2696, 3182)  (2704, 3176)  (2705, 3170)  (2707, 3168)  (2711, 3162)  (2712, 3159)  (2712, 3159)  (2717, 3155)  (2722, 3151)  (2723, 3149)  (2727, 3145)  (2729, 3144)  (2735, 3143)  (2737, 3137)  (2738, 3135)  (2740, 3133)  (2741, 3130)  (2742, 3129)  (2748, 3128)  (2751, 3119)  (2755, 3117)  (2756, 3110)  (2762, 3109)  (2766, 3102)  (2769, 3099)  (2773, 3095)  (2777, 3092)  (2778, 3086)  (2778, 3086)  (2782, 3084)  (2783, 3080)  (2785, 3075)  (2786, 3074)  (2791, 3073)  (2791, 3073)  (2793, 3068)  (2796, 3067)  (2797, 3065)  (2798, 3059)  (2800, 3057)  (2800, 3057)  (2803, 3052)  (2804, 3050)  (2807, 3049)  (2811, 3047)  (2812, 3042)  (2818, 3039)  (2819, 3037)  (2822, 3032)  (2825, 3029)  (2827, 3023)  (2831, 3021)  (2833, 3018)  (2834, 3013)  (2837, 3008)  (2841, 3005)  (2842, 3002)  (2845, 3001)  (2627, 3239)  (2630, 3236)  (2635, 3235)  (2643, 3226)  (2645, 3224)  (2647, 3223)  (2654, 3222)  (2655, 3213)  (2660, 3210)  (2662, 3207)  (2670, 3206)  (2674, 3200)  (2680, 3196)  (2690, 3191)  (2696, 3182)  (2704, 3176)  (2705, 3170)  (2707, 3168)  (2711, 3162)  (2712, 3159)  (2712, 3159)  (2717, 3155)  (2722, 3151)  (2723, 3149)  (2727, 3145)  (2729, 3144)  (2735, 3143)  (2737, 3137)  (2738, 3135)  (2740, 3133)  (2741, 3130)  (2742, 3129)  (2748, 3128)  (2751, 3119)  (2755, 3117)  (2756, 3110)  (2762, 3109)  (2766, 3102)  (2769, 3099)  (2773, 3095)  (2777, 3092)  (2778, 3086)  (2778, 3086)  (2782, 3084)  (2783, 3080)  (2785, 3075)  (2786, 3074)  (2791, 3073)  (2791, 3073)  (2793, 3068)  (2796, 3067)  (2797, 3065)  (2798, 3059)  (2800, 3057)  (2800, 3057)  (2803, 3052)  (2804, 3050)  (2807, 3049)  (2811, 3047)  (2812, 3042)  (2818, 3039)  (2819, 3037)  (2822, 3032)  (2825, 3029)  (2827, 3023)  (2831, 3021)  (2833, 3018)  (2834, 3013)  (2837, 3008)  (2841, 3005)  (2842, 3002)  (2845, 3001)  (2627, 3239)  (2630, 3236)  (2635, 3235)  (2643, 3226)  (2645, 3224)  (2647, 3223)  (2654, 3222)  (2655, 3213)  (2660, 3210)  (2662, 3207)  (2670, 3206)  (2674, 3200)  (2680, 3196)  (2690, 3191)  (2696, 3182)  (2704, 3176)  (2705, 3170)  (2707, 3168)  (2711, 3162)  (2712, 3159)  (2712, 3159)  (2717, 3155)  (2722, 3151)  (2723, 3149)  (2727, 3145)  (2729, 3144)  (2735, 3143)  (2737, 3137)  (2738, 3135)  (2740, 3133)  (2741, 3130)  (2742, 3129)  (2748, 3128)  (2751, 3119)  (2755, 3117)  (2756, 3110)  (2762, 3109)  (2766, 3102)  (2769, 3099)  (2773, 3095)  (2777, 3092)  (2778, 3086)  (2778, 3086)  (2782, 3084)  (2783, 3080)  (2785, 3075)  (2786, 3074)  (2791, 3073)  (2791, 3073)  (2793, 3068)  (2796, 3067)  (2797, 3065)  (2798, 3059)  (2800, 3057)  (2800, 3057)  (2803, 3052)  (2804, 3050)  (2807, 3049)  (2811, 3047)  (2812, 3042)  (2818, 3039)  (2819, 3037)  (2822, 3032)  (2825, 3029)  (2827, 3023)  (2831, 3021)  (2833, 3018)  (2834, 3013)  (2837, 3008)  (2841, 3005)  (2842, 3002)  (2845, 3001)  (2627, 3239)  (2630, 3236)  (2635, 3235)  (2643, 3226)  (2645, 3224)  (2647, 3223)  (2654, 3222)  (2655, 3213)  (2660, 3210)  (2662, 3207)  (2670, 3206)  (2674, 3200)  (2680, 3196)  (2690, 3191)  (2696, 3182)  (2704, 3176)  (2705, 3170)  (2707, 3168)  (2711, 3162)  (2712, 3159)  (2712, 3159)  (2717, 3155)  (2722, 3151)  (2723, 3149)  (2727, 3145)  (2729, 3144)  (2735, 3143)  (2737, 3137)  (2738, 3135)  (2740, 3133)  (2741, 3130)  (2742, 3129)  (2748, 3128)  (2751, 3119)  (2755, 3117)  (2756, 3110)  (2762, 3109)  (2766, 3102)  (2769, 3099)  (2773, 3095)  (2777, 3092)  (2778, 3086)  (2778, 3086)  (2782, 3084)  (2783, 3080)  (2785, 3075)  (2786, 3074)  (2791, 3073)  (2791, 3073)  (2793, 3068)  (2796, 3067)  (2797, 3065)  (2798, 3059)  (2800, 3057)  (2800, 3057)  (2803, 3052)  (2804, 3050)  (2807, 3049)  (2811, 3047)  (2812, 3042)  (2818, 3039)  (2819, 3037)  (2822, 3032)  (2825, 3029)  (2827, 3023)  (2831, 3021)  (2833, 3018)  (2834, 3013)  (2837, 3008)  (2841, 3005)  (2842, 3002)  (2845, 3001)  (2627, 3239)  (2630, 3236)  (2635, 3235)  (2643, 3226)  (2645, 3224)  (2647, 3223)  (2654, 3222)  (2655, 3213)  (2660, 3210)  (2662, 3207)  (2670, 3206)  (2674, 3200)  (2680, 3196)  (2690, 3191)  (2696, 3182)  (2704, 3176)  (2705, 3170)  (2707, 3168)  (2711, 3162)  (2712, 3159)  (2712, 3159)  (2717, 3155)  (2722, 3151)  (2723, 3149)  (2727, 3145)  (2729, 3144)  (2735, 3143)  (2737, 3137)  (2738, 3135)  (2740, 3133)  (2741, 3130)  (2742, 3129)  (2748, 3128)  (2751, 3119)  (2755, 3117)  (2756, 3110)  (2762, 3109)  (2766, 3102)  (2769, 3099)  (2773, 3095)  (2777, 3092)  (2778, 3086)  (2778, 3086)  (2782, 3084)  (2783, 3080)  (2785, 3075)  (2786, 3074)  (2791, 3073)  (2791, 3073)  (2793, 3068)  (2796, 3067)  (2797, 3065)  (2798, 3059)  (2800, 3057)  (2800, 3057)  (2803, 3052)  (2804, 3050)  (2807, 3049)  (2811, 3047)  (2812, 3042)  (2818, 3039)  (2819, 3037)  (2822, 3032)  (2825, 3029)  (2827, 3023)  (2831, 3021)  (2833, 3018)  (2834, 3013)  (2837, 3008)  (2841, 3005)  (2842, 3002)  (2845, 3001)  (2627, 3239)  (2630, 3236)  (2635, 3235)  (2643, 3226)  (2645, 3224)  (2647, 3223)  (2654, 3222)  (2655, 3213)  (2660, 3210)  (2662, 3207)  (2670, 3206)  (2674, 3200)  (2680, 3196)  (2690, 3191)  (2696, 3182)  (2704, 3176)  (2705, 3170)  (2707, 3168)  (2711, 3162)  (2712, 3159)  (2712, 3159)  (2717, 3155)  (2722, 3151)  (2723, 3149)  (2727, 3145)  (2729, 3144)  (2735, 3143)  (2737, 3137)  (2738, 3135)  (2740, 3133)  (2741, 3130)  (2742, 3129)  (2748, 3128)  (2751, 3119)  (2755, 3117)  (2756, 3110)  (2762, 3109)  (2766, 3102)  (2769, 3099)  (2773, 3095)  (2777, 3092)  (2778, 3086)  (2778, 3086)  (2782, 3084)  (2783, 3080)  (2785, 3075)  (2786, 3074)  (2791, 3073)  (2791, 3073)  (2793, 3068)  (2796, 3067)  (2797, 3065)  (2798, 3059)  (2800, 3057)  (2800, 3057)  (2803, 3052)  (2804, 3050)  (2807, 3049)  (2811, 3047)  (2812, 3042)  (2818, 3039)  (2819, 3037)  (2822, 3032)  (2825, 3029)  (2827, 3023)  (2831, 3021)  (2833, 3018)  (2834, 3013)  (2837, 3008)  (2841, 3005)  (2842, 3002)  (2845, 3001)  (2627, 3239)  (2630, 3236)  (2635, 3235)  (2643, 3226)  (2645, 3224)  (2647, 3223)  (2654, 3222)  (2655, 3213)  (2660, 3210)  (2662, 3207)  (2670, 3206)  (2674, 3200)  (2680, 3196)  (2690, 3191)  (2696, 3182)  (2704, 3176)  (2705, 3170)  (2707, 3168)  (2711, 3162)  (2712, 3159)  (2712, 3159)  (2717, 3155)  (2722, 3151)  (2723, 3149)  (2727, 3145)  (2729, 3144)  (2735, 3143)  (2737, 3137)  (2738, 3135)  (2740, 3133)  (2741, 3130)  (2742, 3129)  (2748, 3128)  (2751, 3119)  (2755, 3117)  (2756, 3110)  (2762, 3109)  (2766, 3102)  (2769, 3099)  (2773, 3095)  (2777, 3092)  (2778, 3086)  (2778, 3086)  (2782, 3084)  (2783, 3080)  (2785, 3075)  (2786, 3074)  (2791, 3073)  (2791, 3073)  (2793, 3068)  (2796, 3067)  (2797, 3065)  (2798, 3059)  (2800, 3057)  (2800, 3057)  (2803, 3052)  (2804, 3050)  (2807, 3049)  (2811, 3047)  (2812, 3042)  (2818, 3039)  (2819, 3037)  (2822, 3032)  (2825, 3029)  (2827, 3023)  (2831, 3021)  (2833, 3018)  (2834, 3013)  (2837, 3008)  (2841, 3005)  (2842, 3002)  (2845, 3001)  (2627, 3239)  (2630, 3236)  (2635, 3235)  (2643, 3226)  (2645, 3224)  (2647, 3223)  (2654, 3222)  (2655, 3213)  (2660, 3210)  (2662, 3207)  (2670, 3206)  (2674, 3200)  (2680, 3196)  (2690, 3191)  (2696, 3182)  (2704, 3176)  (2705, 3170)  (2707, 3168)  (2711, 3162)  (2712, 3159)  (2712, 3159)  (2717, 3155)  (2722, 3151)  (2723, 3149)  (2727, 3145)  (2729, 3144)  (2735, 3143)  (2737, 3137)  (2738, 3135)  (2740, 3133)  (2741, 3130)  (2742, 3129)  (2748, 3128)  (2751, 3119)  (2755, 3117)  (2756, 3110)  (2762, 3109)  (2766, 3102)  (2769, 3099)  (2773, 3095)  (2777, 3092)  (2778, 3086)  (2778, 3086)  (2782, 3084)  (2783, 3080)  (2785, 3075)  (2786, 3074)  (2791, 3073)  (2791, 3073)  (2793, 3068)  (2796, 3067)  (2797, 3065)  (2798, 3059)  (2800, 3057)  (2800, 3057)  (2803, 3052)  (2804, 3050)  (2807, 3049)  (2811, 3047)  (2812, 3042)  (2818, 3039)  (2819, 3037)  (2822, 3032)  (2825, 3029)  (2827, 3023)  (2831, 3021)  (2833, 3018)  (2834, 3013)  (2837, 3008)  (2841, 3005)  (2842, 3002)  (2845, 3001)  (2627, 3239)  (2630, 3236)  (2635, 3235)  (2643, 3226)  (2645, 3224)  (2647, 3223)  (2654, 3222)  (2655, 3213)  (2660, 3210)  (2662, 3207)  (2670, 3206)  (2674, 3200)  (2680, 3196)  (2690, 3191)  (2696, 3182)  (2704, 3176)  (2705, 3170)  (2707, 3168)  (2711, 3162)  (2712, 3159)  (2712, 3159)  (2717, 3155)  (2722, 3151)  (2723, 3149)  (2727, 3145)  (2729, 3144)  (2735, 3143)  (2737, 3137)  (2738, 3135)  (2740, 3133)  (2741, 3130)  (2742, 3129)  (2748, 3128)  (2751, 3119)  (2755, 3117)  (2756, 3110)  (2762, 3109)  (2766, 3102)  (2769, 3099)  (2773, 3095)  (2777, 3092)  (2778, 3086)  (2778, 3086)  (2782, 3084)  (2783, 3080)  (2785, 3075)  (2786, 3074)  (2791, 3073)  (2791, 3073)  (2793, 3068)  (2796, 3067)  (2797, 3065)  (2798, 3059)  (2800, 3057)  (2800, 3057)  (2803, 3052)  (2804, 3050)  (2807, 3049)  (2811, 3047)  (2812, 3042)  (2818, 3039)  (2819, 3037)  (2822, 3032)  (2825, 3029)  (2827, 3023)  (2831, 3021)  (2833, 3018)  (2834, 3013)  (2837, 3008)  (2841, 3005)  (2842, 3002)  (2845, 3001)  };
    \addlegendentry{\ {\scriptsize Shared among} {\scriptsize $MOEA_1$ and $MOEA_2$}}
    
 \addplot[samples=4,
         scatter/use mapped color={ball color=blue},
         scatter,
         only marks,
         mark=ball,
         mark size=2.5pt]
   coordinates {   
    (2951,3344)
    };
    \addlegendentry{\ {\scriptsize Ideal point}}

\draw [red!5, fill, opacity = 0.5]
                (axis cs:2842,0) rectangle (axis cs:3040, 3440); 
\draw [yellow!5, fill, opacity = 0.5]
                (axis cs:0,3239) rectangle (axis cs:3040, 3440);  
\draw [gray!20, fill, opacity = 0.5]
                (axis cs:2627,3002) rectangle (axis cs:3040, 3440);

    \draw[dashed, black!80] (axis cs:2951,3344) -- node[left]{} (axis cs:2627, 3239);
    \draw[dashed, black!80] (axis cs:2951,3344) -- node[left]{} (axis cs:2842, 3002);
    \draw[dashed, black!80] (axis cs:2951,3344) -- node[left]{} (axis cs:2277, 3344);
    \draw[dashed, black!80] (axis cs:2951,3344) -- node[left]{} (axis cs:2951, 2651);
    \draw[ gray!30] (axis cs:2842,0) -- node[left]{} (axis cs:2842, 3460);
    \draw[ gray!30] (axis cs:2627,0) -- node[left]{} (axis cs:2627, 3460);
	\draw[ gray!30] (axis cs:0,3239) -- node[left]{} (axis cs:3040, 3239);
	\draw[ gray!30] (axis cs:0,3002) -- node[left]{} (axis cs:3040, 3002);

 \end{axis}

\end{tikzpicture}

\caption{{\small Illustrative example of the partitioning procedure.}}
 \label{fig: parti}
\end{figure}
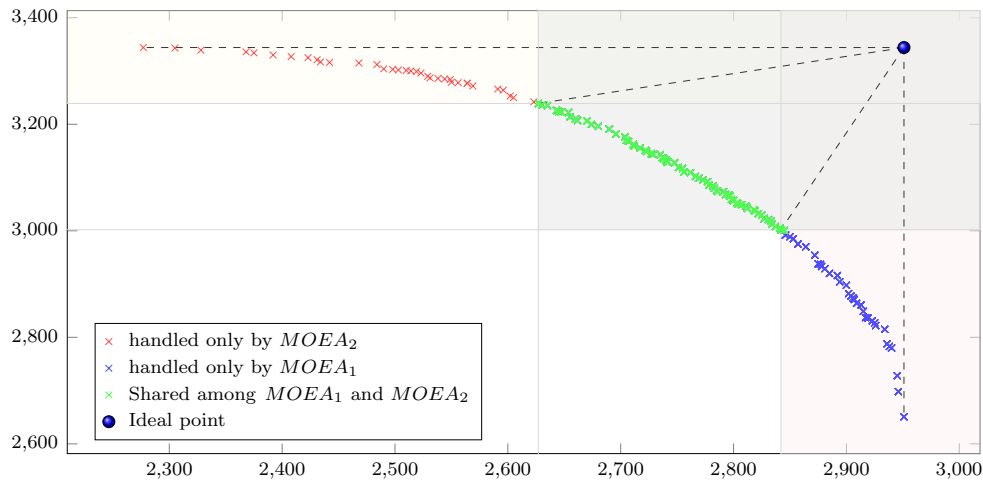

  Figure \ref{fig: parti} presents an example of the partitioning procedure applied to a bi-objective Knapsack instance: 2KP100-50 \cite{9}, showing the obtained Pareto front at iteration $10^3$, partitioned using $\alpha = 0.25$. Here the shared solutions (i.e., contained in the gray area) are enclosed in the interquartile range (IQR) of the Pareto solutions. Furthermore, IQR is represents here $50\%$ of the current non-dominated solutions that are shared among search entities.%

\section{\textbf{Experimental studies}}
We tested the sugessted algorithm on benchmark instances of MOMKP chosen from the instance libraries: Zitzler and al. \cite{16}, of which we consider for this experiments six instances with the number of items 250, 500, and 750, with two and three objectives. The  parallelization of the suggested algorithm is implemented via multi-threading under JAVA SE platform. The tests was carried out on a personal computer equipped with an Intel\textsuperscript{{\tiny \textregistered}} Core\textsuperscript{{\tiny TM}} i7-5600U CPU, 2.60GHz with 8GBs of RAM. We compared the performance of the PCPMOEA with four multiobjective algorithms with different concepts and/or different search strategies: 
\begin{itemize}
\item NSGAII \cite{6}: An elitist non-dominated sorting genetic algorithm, a multiobjective genetic algorithm using dominance depth, and crowding distance for selection operator and search guidance. 
\item SPEA2 \cite{7}: Strength Pareto evolutionary algorithm, a MOEA with an external archive using dominance rank, and a nearest neighbor density estimation technique. 
\item MOEA/D \cite{37}: Multiobjective evolutionary algorithm based on decomposition, essentially, it  decomposes an MOP into a number of scalar optimization subproblems and optimizes them simultaneously.
\item MOFPA \cite{5}: Multiobjective Firefly algorithm with Particle swarm optimization, a hybrid swarm intelligence discrete algorithm, employing cooperation of two intelligent swarm algorithms: Firefly Algorithm, and Particle swarm optimization. 
\end{itemize}
In order to evaluate and compare the quality of solutions (convergence, and the spread) evoloved using these algorithms, we used four performance metrics: Inverted Generational Distance (IGD) \cite{38} , Hypervolume \cite{39}, the Spacing metric \cite{36}, and the set coverage metric \cite{39}. As we mentioned, one of the aims of the suggested algorithm is to minimize the distance between the potentially efficient solutions and the ideal point. Therefore, we used a modified version of the GD metric as to measure this distance. \vspace{-0.5em}
\subsection{Performance evaluation metrics\vspace{-0.5em}}
\subsubsection*{Hypervolume \cite{39}}
One of the most popular indicators for multiobjective optimization algorithms is the hypervolume . The hypervolume indicator measures the volume of the $k$-dimensional space dominated by a set of points $A$ (in the objectives space). The aim of this indicator is to measure both the convergence to the true Pareto front and diversity of the obtained Pareto fronts. The calculation of this volume requires the designation of a reference point $Z_{ref}$, which is dominated by all the points of the set $A$. Consequently, a set with a larger hypervolume is likely to present a better set of trade-offs than sets with lower hypervolume. The reference point considered to compute the hypervolume is the null vector (origin point, $Z_{ref}= \{0\}^{k}$, where $k$ is the number of objective functions). 
\subsubsection*{Inverted generational distance (IGD) \cite{35}}
The IGD metric is an inverted version of the GD (generational distance metric) used to measure convergence. This metric evaluates the average distance between the approximated set $A$ and a reference set PF (true Pareto front). In other words, the IGD measures the proximity between the obtained potentially efficient solutions and the true Pareto front in the objective space.  
$$\displaystyle IGD =  \frac{\left( \sum_{s\in PF} d^2_s \right )^{\frac{1}{2}}}{|PF|},$$
where, $d_s$ is the Euclidean distance in the objective space between the solution $s\in PF$ and the nearest solution in $A$.
$$\forall s \in PF, \; d_s = \min_{ s'\in A }{|| Z(s) - Z(s') ||}_{2}.$$
\subsubsection*{Ideal distance (ID)}
As to evaluate the distance between the obtained potentially efficient solutions and the ideal vector, we used a modified version the GD metric, we refer to it as the ideal distance (ID). Here, we compute the average distance between the approximated set $A$ and the ideal vector $Z_0$. Similar to the GD metric, the ID is calculated as follows:
$$\displaystyle ID \overset{\Delta}{=}  \frac{\left( \sum_{s\in A} d^2_s \right )^{\frac{1}{2}}}{|A|},$$
where, $\forall s \in A, \; d_s = {|| Z(s) - Z_0 ||}_{2}.$
\subsubsection*{Spacing metric (SP) \cite{36}}
This metric evaluates the spread of a given set of non-dominated solutions, say $A$, using its distribution in the objective space. This is, by computing the standard deviation in the distances between consecutive pairs of solutions. 
$$SP = \sqrt{\frac{1}{|A|-1} \sum_{s\in A} {(\bar{d}-d_s)}^2},$$
where, $d_s$ represents the distance between the solution $s$ and the closest neighboring solution in the objective space.
$$\forall s \in A, \; d_s = \min_{s' \in A, s'\not=s} \sum_{m=1}^k |Z^m(s)-Z^m(s')|, $$ 
and $\bar{d}$ denotes the average of the distances $d_s,$  $\forall s\in A$.
However, a set $A$ with a good spacing metric does not necessarily mean a good distribution of solutions compared to the entire Pareto optimal front \cite{40}.
\subsubsection*{ Set coverage metric \cite{39}} This metric is used for comparing two potentially efficient solutions sets, say $A$ and $B$. This is, by calculating the proportion of solutions in the set $B$ that are weakly dominated by solutions in the first set $A$. 
$$ C(A,B)=\frac{|\{ b\in B | \exists a\in A : a \succeq b\}|}{|B|}\cdot $$
Note that $C(A, B)$ and $C(B, A)$ are not necessarily equal, hence, both must be calculated in order to compare the approximations $A$ and $B$.
\subsection{Experimental results and comments}
In the rest of this section we will present the experimental results obtained for the following algorithms: PCPMOEA, NSGAII, SPEA2, MOEA/D, and MOFPA. Furthermore, for each algorithm we used fifteen potentially efficient solutions (15 runs for each instance). The aim of these experimentations is to compare the suggested algorithm (PCPMOEA) with a set of efficient algorithms that uses different approaches of selection mechanisms and search strategies. The parameters we used for each algorithm are chosen according to the origin papers, and for the PCPMOEA the parameters are fixed as follows: population size $N = 150$, partitioning parameter $\alpha = 0.25$ for $2$ objectives,  $N = 250$,  $\alpha = 0.3$  for $3$ objectives, with different values for the period. Crossover probability: $0.8$, mutation probability: $0.1$.

As to illustrate the obtained results used for comparison, we qualitatively show in Figures 4-8 some of the obtained results, where each figure shows Pareto fronts obtained for each instance using the above mentioned algorithms. As to maintain the clarity of the three dimensional illustrations, we have chosen to show only the most competitive Pareto fronts, which lead us to the Pareto fronts evolved using PCPMOEA and MOFPA.

Table 1 resumes the obtained results regarding the Hypervolume indicator. The aim of this experimentation is to evaluate and compare the quality of the suggested algorithm with the above mentioned algorithms. The obtained results show that the PCPMOEA algorithm achieves a good quality of solutions, especially, when compared to the NSGAII, SPEA2, and MOEA/D. The PCPMOEA demonstrates to have a better Hypervolume with a significant difference. However, MOFPA appears to be the most competitive to the PCPMOEA, namely, for the instance with three objectives and 750 items (3.750). \vspace{-3mm}

\begin{figure}[H]
\begin{footnotesize}
\begin{center}
\centering
\setlength{\tabcolsep}{.7em}
{\renewcommand{\arraystretch}{1.05}
\begin{tabular}{llrrrrr}
\multirow{ 2}{*}{Instance} &  & \multicolumn{ 5}{c}{Algorithm}  \\\cline{3-7}
 & & NSGAII & SPEA2 & MOEA/D & MOFPA & PCPMOEA \\ \hline
\multirow{ 1}{*}{2.250} 
 & Average &8.76512 & 9.18159 & 9.83537& 9.85685 &\textbf{9.86569}\\
 & Median & 8.77438& 9.18573 & 9.83568 & 9.85572&\textbf{9.86649}\\
 & std &2.716E-2 & 6.837E-2 &6,258E-3& \textbf{2.808E-3}&6.579E-3\\
 & Best & 8.81628 & 9.26775 &9.84387 & 9.86246&\textbf{9.87017}\\
 & Worst & 8.72353 & 9.04351& 9.82125& \textbf{9.85441} & 9.84700 \\
 \hline
 \multirow{ 1}{*}{2.500} 
 & Average & 34.22430 &36.65644  & 40.52875& 40.71232 &\textbf{40.74421}\\
 & Median & 34.24338 &  36.67221& 40.53538& 40.70868& \textbf{40.76286}\\
 & std & 0.1296 & 0.25230 &0.003981 &\textbf{0.00230}&0.00358\\
 & Best & 34.40142 & 36.95143& 40.59179&40.75594& \textbf{40.78454}\\
 & Worst & 34.24338 & 36.19984 & 40.53538 & 40.66021& \textbf{40.69426} \\
 \hline

 \multirow{ 1}{*}{2.750} 
 & Average & 71.42760 & 77.37992 & 88.49387& 88.70273 & \textbf{89.20554}\\
 & Median &71.41220 & 77.33868 & 88.53043& 88.70192 & \textbf{89.22105}\\
 & std & 0.503297 &  0.403669& 9.422E-2 & \textbf{4.897E-2 }& 8.0114E-2\\
 & Best &72.95212 & 78.03857& 88.61037 & 88.81874 &\textbf{89.27410}\\
 & Worst & 70.91351  & 77.33868 & 88.32110 & 88.63514 & \textbf{88.99637}\\
 \hline
 \multirow{ 1}{*}{3.250} 
 & Average & 72322.4 & 77614.4& 92409.8 & 93009.9 & \textbf{93119.6}\\
 & Median &72311.1 & 77443.4 & 92434.8 & 93016.2 & \textbf{93126.2}\\
 & std &435.1785 & 709.9357 &  80.0224& 39.3793& \textbf{37.8173}\\
 & Best & 73187.1 & 79520.4& 92534.5 & 93071.6 & \textbf{93182.2}\\
 & Worst &71477.3 & 76865.1 & 92434.8 & 92956.6 & \textbf{93054.2} \\
 \hline
 \multirow{ 1}{*}{3.500} 
 & Average & 547738.3 & 605712.2 & 758240.9 & 758312.6 & \textbf{766498.7}\\
 & Median &  546821.8 & 605755.9 & 758202.0 & 758331.7 & \textbf{766488.1}\\
 & std & 3608.78 & 3340.01& \textbf{670.63} & 870.71 & 1014.90\\
 & Best &    555481.1 & 612315.4 & 759647.3 & 759959.4& \textbf{768077.1}\\
 & Worst &   542001.8 & 600666.8 & 756955.4 & 756555.1 & \textbf{764012.6}\\
 \hline
 \multirow{ 1}{*}{3.750} 
 & Average & 2195528.6 & 2207060.2 & 2649034.1& 2703273.4 & \textbf{2703522.7}\\
 & Median  & 2192486.3 & 2206634.2 & 2648629.3 & 2703114.3 & \textbf{2703648.3}\\
 & std & 13524.87 & 7704.53& 2214.50 & 732.79 & \textbf{580.84}\\
 & Best    & 2217966.1 & 2225107.6 & 2653975.2 & \textbf{2704278.8} & 2704231.5\\
 & Worst   & 2171692.9 & 2191793.4 & 264862.93 & 2702333.1& \textbf{2702581.8}\\
 \hline
 \multicolumn{7}{l}{Values in bold indicate the best obtained result in the correspondent row.}
\end{tabular}
}
\end{center}
\end{footnotesize}
\capt{1}{Experimental results concerning the Hypervolume indicator of the MOMKP instances (the unit of values is $10^{7}$).}
\label{hyperv}
\end{figure}

In order to assess the convergence of the obtained Pareto fronts towards the true Pareto front, we computed the IGD metric values for all fifteen runs of each algorithm. The obtained results are resumed in Table 2. From Table 2 it is obvious that the obtained average (and median) values of IGD for each instance using PCPMOEA are better than those obtained by the other four algorithms. Accordingly, the PCPMOEA has converged better towards the true Pareto front compared to other tested algorithms. These results can be also confirmed qualitatively by observing Figures 4-8.  

\begin{figure}[H]
\begin{footnotesize}
\begin{center}
\centering
\setlength{\tabcolsep}{.7em}
\setlength{\arrayrulewidth}{.05em}
{\renewcommand{\arraystretch}{1.05}
\begin{tabular}{llrrrrr}
\multirow{ 2}{*}{Instance} &  & \multicolumn{ 5}{c}{Algorithm}  \\\cline{3-7}
 & & NSGAII & SPEA2 & MOEA/D & MOFPA & PCPMOEA \\ \hline
\multirow{ 1}{*}{2.250} 
 & Average & 76.09615 & 15.64012 &4.00665 &1.07286  &\textbf{0.29718}\\
 & Median &  72.33389 & 15.63951 &  3.94590 & 1.05223 &\textbf{0.29637}\\
 & std &9.88158& 2.25782 & 0.39901 &\textbf{0.10124} & 0.12132\\
 & Best & 61.35042 & 12.16928& 3.13348 & 0.89778 & \textbf{0.12869}\\
 & Worst & 95.26060 &20.76262  & 4.86990 & 1.26336 & \textbf{0.57562}\\
 \hline
%
%
%

 \multirow{ 1}{*}{2.500} 
 & Average & 355.87107 & 81.02844 & 14.67452& 2.18277 &\textbf{0.69520}\\
 & Median & 359.74476 & 80.88057& 14.67389 &2.17060 & \textbf{0.58288}\\
 & std & 43.97135 & 12.46397 & 1.01867& \textbf{0.23388} & 0.29223\\
 & Best & 293.13791& 65.77132 & 12.72309 & 1.84681 & \textbf{0.37414}\\
 & Worst &418.81440 & 111.14670 & 16.61654 & 2.64536 & \textbf{1.35093}\\
 \hline
 
%
%

 \multirow{ 1}{*}{2.750} 
 & Average & 616.90264 & 81.028445 & 32.19721 & 10.54516 &\textbf{3.53590}\\
 & Median & 622.15216& 80.88057 & 32.51318& 10.94526 & \textbf{3.33629}\\
 & std & 55.24098& 12.46397 & 1.48806& \textbf{0.62207}& 0.90879\\
 & Best &486.57075 & 65.77132 & 29.99165& 9.70699 & \textbf{2.37712}\\
 & Worst & 722.35960& 111.14670 & 34.86158& 11.02817& \textbf{6.22080}\\
 \hline
 \multirow{ 1}{*}{3.250} 
 & Average &45.52375 & 16.23703 & 6.4468998 & 2.3051193 & \textbf{1.3234943}\\
 & Median &45.82247 & 16.33518 & 6.4847662 &  2.3301823& \textbf{1.2882777}\\
 & std & 3.11597 & 1.36025 & 0.2016878&\textbf{0.0500886}& 0.1688584\\
 & Best & 39.93595 & 13.60152 &  6.0662901& 2.2376310& \textbf{1.0267159}\\
 & Worst & 51.34482& 18.37225  &6.7613303 &  2.3671946& \textbf{1.6474703}\\
 \hline
%
%
%
%

 \multirow{ 1}{*}{3.500} 
 & Average & 159.08672 & 47.10854 & 14.4351057 & 4.7504078 & \textbf{3.7412855}\\
 & Median  & 158.70766 & 47.61122 & 14.4262958 & 4.8048472 & \textbf{3.3671623}\\
 & std     & 11.05657  & 4.19196 &  0.5505419 & \textbf{0.1500751} & 1.1078586\\
 & Best    & 135.34139 & 37.73514 & 13.6882246 & 4.4624758 & \textbf{2.4686905}\\
 & Worst   & 174.41828 & 54.96809 & 15.6030902 & \textbf{4.9110414} & 6.9519051\\
 \hline
%
%
%

 \multirow{ 1}{*}{3.750} 
 & Average & 59.62771 & 55.42584 &  26.7652333 &3.9739178 &\textbf{2.6863341}\\
 & Median &  59.82934&  55.06092 & 26.6458393 & 3.9824589 & \textbf{2.6076377}\\
 & std & 3.16625& 3.49813 &  0.5282267 &\textbf{ 0.0337948 }& 0.2612799\\
 & Best &    52.14942& 47.98128 & 25.9533346 & 3.9126965 & \textbf{2.4197297}\\
 & Worst & 64.00565 &  61.57806 & 27.9284213 & 4.0265507 & \textbf{3.2129786} \\
 \hline
 \multicolumn{7}{l}{Values in bold indicate the best obtained result among the correspondent row.}
\end{tabular}
}
\end{center}
\end{footnotesize}
\capt{2}{Experimental results concerning the IGD of the MOMKP instances.}
\label{igd}
\end{figure}

Table 3 summarizes the results regarding the spacing metric values, obtained for all five algorithms using all fifteen runs for each algorithm. The aim of this experimentation is to compare the diversity (i.e., the distribution uniformity) of the produced Pareto fronts. The results show that all of the competing algorithms can produce Pareto fronts with a good spacing. However, the Pareto fronts evolved by the proposed algorithm have shown to be the most uniformly distributed for almost all the tested instances. Namely for the following instances: 2.500, 2.750, 3.250, 3.500, and 3.750. Consequently, we can say that the diversity of Pareto fronts evolved by PCPMOEA are at least comparable to those of the other algorithms.

\begin{figure}[H]
\begin{footnotesize}
\begin{center}
\centering
\setlength{\tabcolsep}{.7em}
\setlength{\arrayrulewidth}{.05em}
{\renewcommand{\arraystretch}{1.05}
\begin{tabular}{llccccc}
\multirow{ 2}{*}{Instance} &  & \multicolumn{ 5}{c}{Algorithm}  \\\cline{3-7}
 & & NSGAII & SPEA2 & MOEA/D & MOFPA & PCPMOEA \\ \hline
\multirow{ 1}{*}{2.250} 
 & Average & 4.8441E-3& 1.6520E-3& 1.6788E-3& 1.8034E-3 & 1.0830E-3\\
 & Median & 4.1695E-3& 1.4919E-3& 1.7452E-3&  1.7452E-3 & 1.0712E-3\\
 & std & 2.7102E-3& 5.5580E-4& 2.9105E-4& 2.9105E-4 & 3.9267E-4\\
 & Best & 2.0031E-3& 1.2372E-3& 1.1081E-3& 9.7163E-4 & 4.3610E-4\\
 & Worst & 1.1410E-2& 3.2591E-3& 2.0978E-3& 2.6706E-3& 1.7888E-3\\
 \hline

 \multirow{ 1}{*}{2.500} 
 & Average & 5.3383E-3&1.5020E-3 & 1.0604E-3& 9.9709E-4& 4.5168E-4\\
 & Median & 5.4688E-3& 1.4126E-3&  9.9696E-4& 9.6733E-4&4.2077E-4\\
 & std & 3.2216E-3& 4.3462E-4& 1.5039E-4& 1.5387E-4&  1.1816E-4\\
 & Best & 1.5329E-3& 9.4383E-4& 9.1382E-4& 7.8784E-4& 2.5867E-4\\
 & Worst & 1.2600E-2& 2.6173E-3& 1.4491E-3& 1.4351E-3& 6.6605E-4\\
 \hline
 \multirow{ 1}{*}{2.750} 
 & Average & 3.5060E-3& 1.5047E-3& 1.0523E-3& 1.4097E-3& 4.4397E-4\\
 & Median & 2.8953E-3& 1.4157E-3& 1.0648E-3& 1.3104E-3& 4.3931E-4\\
 & std & 1.4508E-3& 4.1546E-4& 6.5984E-5& 3.6208E-4& 7.3312E-5\\
 & Best & 1.8880E-3& 9.6539E-4& 8.9261E-4& 9.0932E-4& 3.3732E-4\\
 & Worst & 6.6667E-3& 2.1664E-3& 1.1373E-3& 1.8850E-3& 6.4364E-4\\
 \hline
%
 \multirow{ 1}{*}{3.250} 
  & Average & 4.5967E-3& 2.5863E-3& 2.6694E-3& 2.1473E-3& 1.2423E-3\\
 & Median & 4.5748E-3& 2.5642E-3& 2.6888E-3 & 2.1386E-3& 1.2427E-3\\
 & std & 8.4339E-4& 3.6247E-4& 8.1093E-5 & 1.0228E-4& 6.9140E-5\\
 & Best & 3.5991E-3  & 1.9196E-3& 2.5017E-3& 1.9909E-3& 1.9909E-3\\
 & Worst & 6.3322E-3& 3.1437E-3& 2.8073E-3 & 2.3278E-3& 2.3278E-3\\
 \hline
 
 \multirow{ 1}{*}{3.500} 
 & Average & 3.8769E-3& 1.9653E-3& 2.3068E-3 & 2.2222E-3 & 1.2039E-3 \\
 & Median &  3.7883E-3& 1.9151E-3& 2.3036E-3& 2.2210E-3& 1.1587E-3\\
 & std & 6.9103E-4& 2.4949E-4& 3.8202E-5 & 8.0886E-5& 1.7005E-4\\
 & Best & 2.7061E-3& 1.5828E-3&2.2471E-3 & 2.1169E-3& 1.0065E-3\\
 & Worst & 5.5471E-3&2.3522E-3&2.3762E-3 & 2.4065E-3& 1.6452E-3\\
 \hline
%
 \multirow{ 1}{*}{3.750} 
 & Average & 1.6451E-3& 9.7395E-4& 2.2277E-3 & 1.2556E-3 & 7.1063E-4\\
 & Median & 1.6529E-3& 9.9614E-4 & 2.2156E-3 &1.2700E-3 & 7.0007E-4\\
 & std &8.5805E-5& 1.0503E-4 & 7.3563E-5 & 7.3563E-5 & 2.6296E-5\\
 & Best & 1.5284E-3& 6.9086E-4 & 2.0839E-3 & 2.0839E-3 & 6.8251E-4\\
 & Worst & 1.7873E-3 & 1.1137E-3& 2.3914E-3 & 2.3914E-3& 7.4769E-4\\
 \hline
\end{tabular}
}
\end{center}
\end{footnotesize}
\capt{3}{Experimental results concerning the spacing metric of the MOMKP instances.}
\label{sp}
\end{figure}

Aspiring to achieve a good convergence towards the true Pareto front, the suggested algorithm operates with multiple MOEAs, within each one, a selection operator is designed to take into consideration one sole criterion through the search process. This is, as we mentioned above, to minimize the distance between the current Pareto front and the ideal point. However, as to evaluate this distance we defined a modified version of the GD metric, called the ideal distance (ID). Table 4 summarizes the obtained results regarding the ID metric. The obtained mean values clearly shows that the Pareto fronts evolved using the proposed algorithm are the nearest to the ideal point.

\begin{figure}[H]
\begin{footnotesize}
\begin{center}
\centering
\setlength{\arrayrulewidth}{.05em}
{\renewcommand{\arraystretch}{1.05}
\begin{tabular}{llrrrrr}
\multirow{ 2}{*}{Instance} &  & \multicolumn{ 5}{c}{Algorithm}  \\\cline{3-7}
 & & NSGAII & SPEA2 & MOEA/D & MOFPA & PCPMOEA \\ \hline

 
\multirow{ 1}{*}{2.250} 
 & Average &495.30 & 349.40 & 213.24 & 206.26 &\textbf{137.97}\\
 & Median & 344.28 &  344.28 &  214.67 & 206.26 &\textbf{136.11}\\
 & std &11.54 & 11.54 &9.18 &   \textbf{4.89} &6.79\\
 & Best & 333.77 & 333.77 &   197.83 &   200.19 & \textbf{129.23}\\
 & Worst & 369.89 & 369.89&    226.26&   217.28 & \textbf{152.99} \\
 \hline
 

 \multirow{ 1}{*}{2.500} 
 & Average & 795.7 & 669.99  & 336.78 & 300.33 &\textbf{190.83}\\
 & Median & 810.5 &  668.75 & 338.80 & 299.53 & \textbf{184.77}\\
 & std        & 132.0 & 68.69 &  8.94 & \textbf{7.60} & 17.68\\
 & Best & 588.5 & 567.98& 567.98 & 288.13 & \textbf{165.24}\\
 & Worst & 1003.4 & 842.07 & 349.65 & 316.11 &  \textbf{229.60}\\
 \hline


 \multirow{ 1}{*}{2.750} 
 & Average & 1145.6 & 1071.2 & 525.94 & 464.08 & \textbf{328.33} \\
 & Median & 1111.7 & 1063.0 & 526.68 & 465.33 & \textbf{329.61}\\
 & std & 12.24 &  10.13 & 10.61 & \textbf{11.94} & 32.88\\
 & Best & 908.1 & 952.5 & 507.04 & 440.59 & \textbf{286.80}\\
 & Worst & 1325.3  & 1325.3 & 542.69 & 486.11 & \textbf{419.83}\\
 \hline
 
 
 \multirow{ 1}{*}{3.250} 
 & Average & 208.35 & 135.82 & 70.75 & 71.43 & \textbf{42.20}\\
 & Median & 213.46 & 135.53 & 70.84 & 71.27 & \textbf{42.36}\\
 & std & 14.80 & 8.54 &  1.48 & \textbf{0.63} & 2.12\\
 & Best & 183.87 & 124.22 & 68.35 & 69.58 & \textbf{38.69}\\
 & Worst & 231.83 & 154.58 & 74.20 & 72.3 & \textbf{46.72} \\
 \hline


 \multirow{ 1}{*}{3.500} 
 & Average & 378.24 & 230.40 & 110.95 & 120.60 & \textbf{86.50}\\
 & Median &  377.71 & 229.35 & 110.77 & 121.18 & \textbf{84.11}\\
 & std & 25.53 & 8.51 & 1.62 & \textbf{2.61} & 11.06\\
 & Best &    335.10 & 217.55 & 108.23 & 115.95 & \textbf{71.81}\\
 & Worst &   434.09 & 245.55 & 113.60 & 124.07 & \textbf{113.87}\\
 \hline


 \multirow{ 1}{*}{3.750} 
 & Average & 413.47 & 417.46 & 161.56 & 117.38 & \textbf{86.56}\\
 & Median  & 417.85 & 417.85 & 161.79 & 117.03 & \textbf{86.56}\\
 & std & 3.16 & 3.16 & 3.25 & \textbf{0.90} & 4.36\\
 & Best    & 411.34 & 411.34 & 156.38 & 116.55 & \textbf{81.70}\\
 & Worst   & 422.11 & 422.11 & 166.41 & 118.88 & \textbf{94.73}\\
 \hline
 \multicolumn{7}{l}{ Values in bold indicate the best obtained result in the correspondent row.}
\vspace{-2mm}
\end{tabular}
}
\end{center}
\end{footnotesize}
\capt{4}{Experimental results concerning the ideal distance of the MOMKP instances.}
\label{id}
\vspace{-2mm}
\end{figure}
As to confirm and support the accuracy of the obtained results regarding the comparison of the convergence and diversity of Pareto fronts, we computed the coverage of each pair of Pareto fronts: (PCPMOEA, and a competing algorithm), produced by all 15 runs of each algorithm. Table 5 shows the obtained mean coverage values for each pair adduced as follows: the symboles $\succeq$ and $\preceq$ refer to C(PCPMOEA, Algorithm) and C(Algorithm, PCPMOEA) respectively. The results show that PCPMOEA produces a better quality of Pareto fronts when compared to NSGAII, SPEA2, MOEA/D. However, MOFPA is shown to be the most competitive to the suggested algorithm, especially for instances: 2.750, 3.250, 3.500, 3.750, although, the suggested algorithm maintained to be dominant over MOFPA, scoring an overall mean coverage values of $52.7\%$ for PCPMOEA against $25.9\%$ for MOFPA.   
\begin{figure}
\centering
\setlength{\tabcolsep}{0.3em}
\setlength{\arrayrulewidth}{.05em}
{\renewcommand{\arraystretch}{1.2}
\begin{footnotesize}
\begin{tabular}{cccccc}
\multirow{ 2}{*}{Instance} &  & \multicolumn{ 4}{c}{Algorithm}  \\\cline{3-6}
& & NSGAII & SPEA2 & MOEA/D & MOFPA \\ \hline
\multirow{ 1}{*}{2.250} & $\succeq$ & 1& 1& 0.9973 &  0.7421\\
                        & $\preceq$ & 0& 0& 0.0003&  0.0819\\ \hline
\multirow{ 1}{*}{2.500} & $\succeq$ & 1& 1& 0.9997 & 0.5859\\
                        & $\preceq$ & 0& 0& 0& 0.2013\\ \hline
\multirow{ 1}{*}{2.750} & $\succeq$ & 1& 1& 1& 0.3834\\
                        & $\preceq$ & 0& 0& 0& 0.2092\\ \hline
\multirow{ 1}{*}{3.250} & $\succeq$ & 1& 1& 0.9278& 0.4525 \\
                        & $\preceq$ & 0& 0& 2.7E-5& 0.2761\\ \hline
\multirow{ 1}{*}{3.500} & $\succeq$ & 1& 1&0.9563 & 0.5671\\
                        & $\preceq$ & 0& 0& 0& 0.4028\\ \hline
\multirow{ 1}{*}{3.750} & $\succeq$ & 1& 1& 1& 0.4332 \\
                        & $\preceq$ & 0& 0& 0& 0.3840\\ \hline\hline                                             
\multirow{ 2}{*}{Average} & $\succeq$ & 1 & 1 & 0.9801& 0.5273\\
						& $\preceq$ & 0 &	0 & 5.45E-5 & 0.2592\\ \hline
\multicolumn{6}{l}{{\scriptsize The last two rows contain the mean values for each column. }}
\end{tabular}
\end{footnotesize}
}
\capt{5}{ Coverage matric of PCPMOEA against other algorithms.}
\label{cov}
\end{figure}

As it is mentioned earlier, in order to visually compare the quality of non-dominated solutions obtained using the five algorithms: NSGAII, SPEA2, MOEA/D, MOFPA, and PCPMOEA. The obtained Pareto fronts for the 6 tested instances are shown in Figures 4-8. From these observations we can confirm the test results regarding the convergence and the spread of the Pareto fronts evolved using PCPMOEA, which can be seen for instances 2.500, 2.750, 2.250, 3.500. It also confirms the fact that MOFPA is the most competitive to the suggested algorithm. Furthermore, we can also observe especially in Figure 8 (instance 3.750), that the suggested algorithm tends to neglect some subspaces with only two conflicting criteria. This is because the chosen selection mechanisms for PCPMOEA tends to handle either one (for each parallel MOEA) or all objective functions combined using dominance relation (for each parallel MOEA, and for the master's selection operator).

\begin{figure}[H]
\includegraphics[scale=0.36]{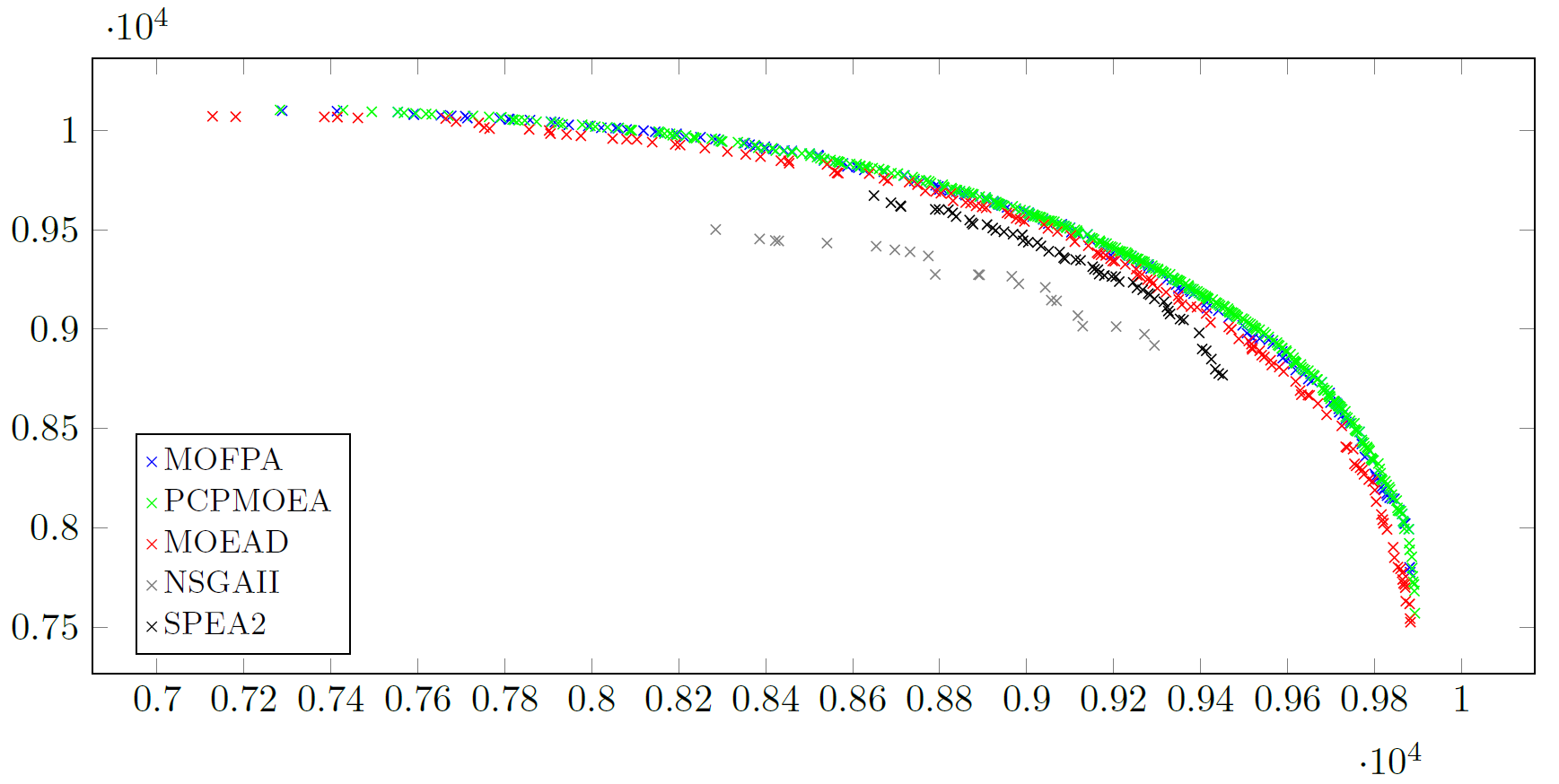}\vspace{-2mm}
\caption{{\footnotesize Illustion of Non-dominated solutions obtained for 2.250 instance.}}
\end{figure}
\begin{figure}[H]
\includegraphics[scale=0.35]{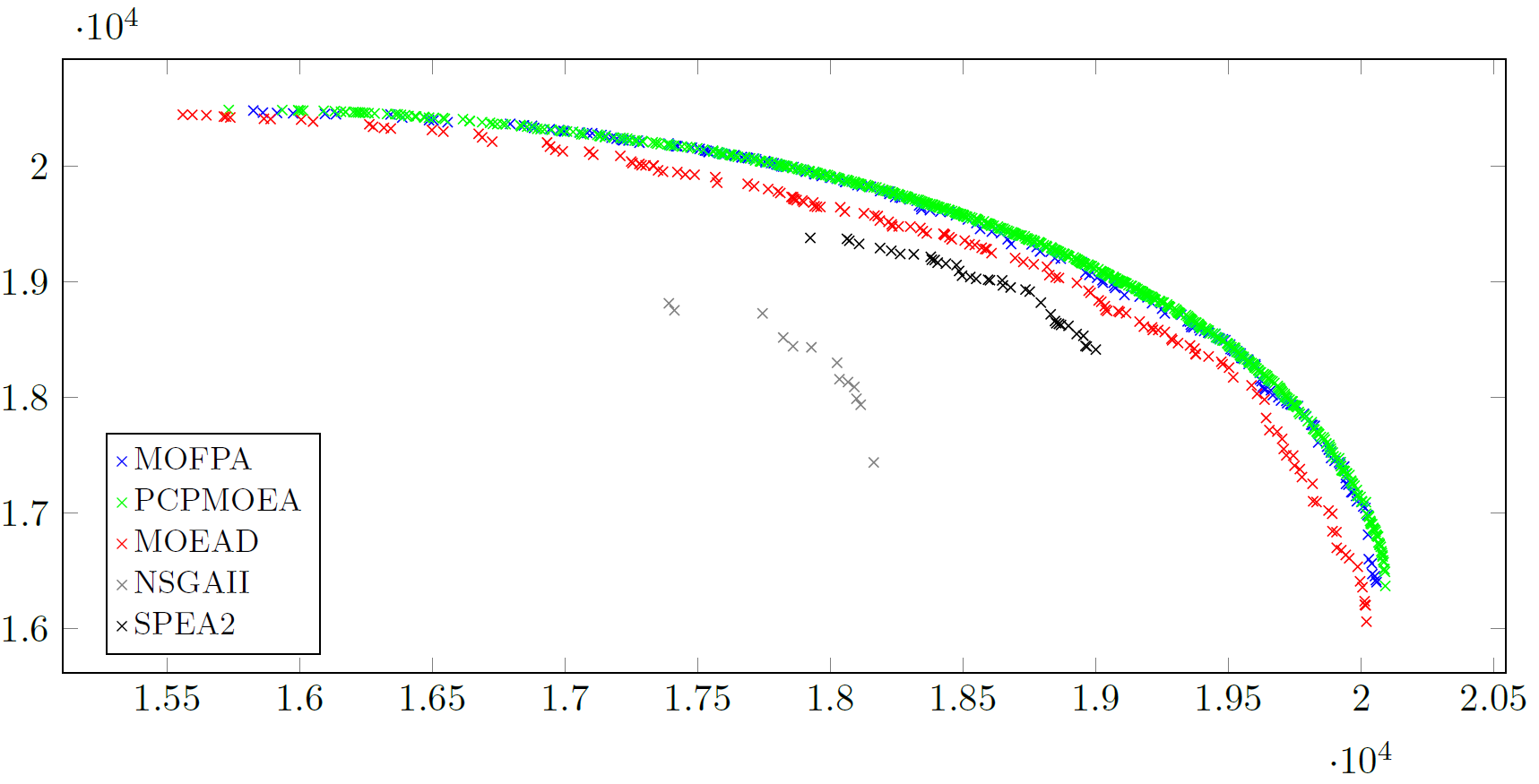}\vspace{-3mm}
\caption{{\footnotesize Illustion of Non-dominated solutions obtained for 2.500 instance.}}
\end{figure}
\begin{figure}[H]
\includegraphics[scale=0.35]{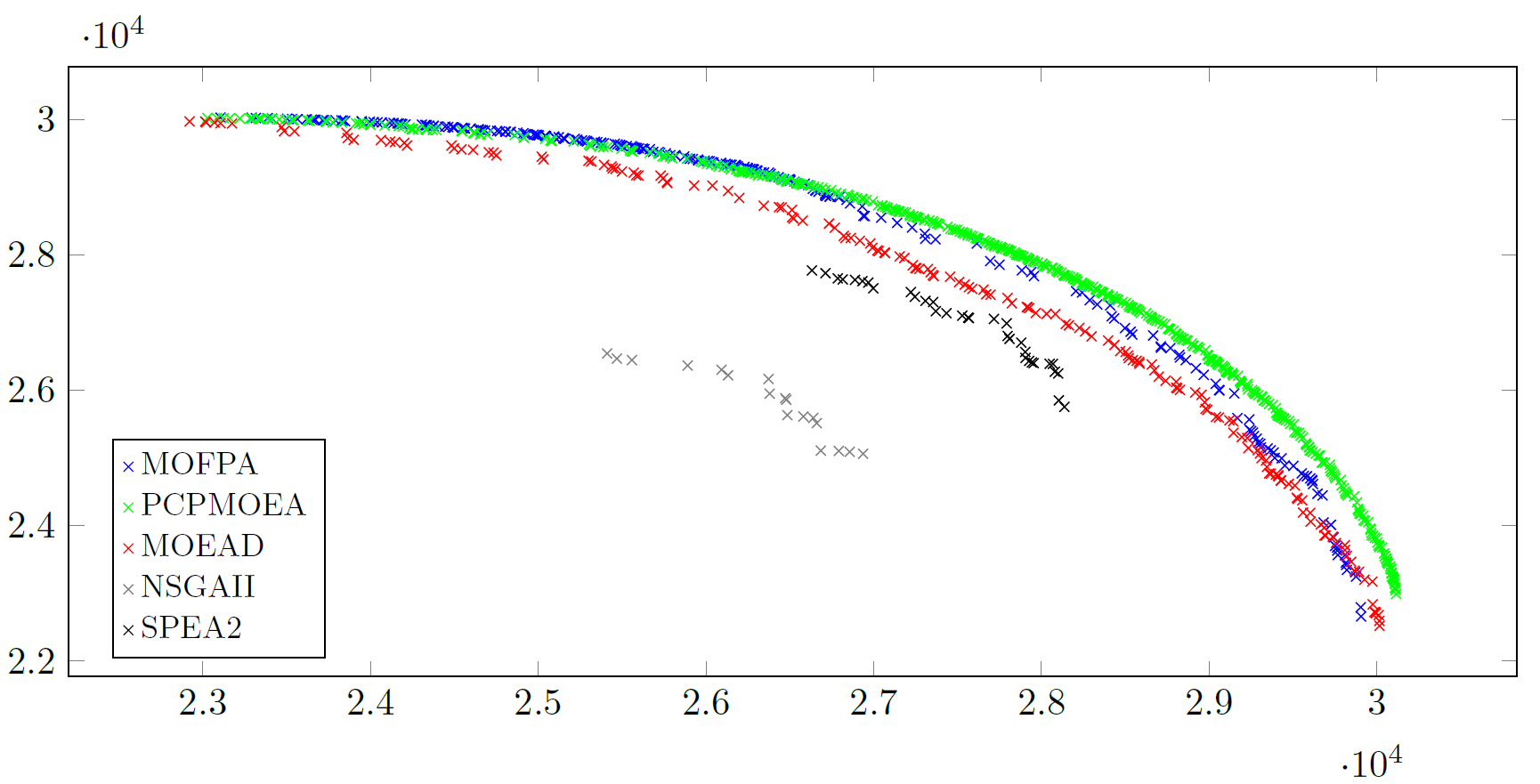}\vspace{-3mm}
\caption{{\footnotesize Illustion of Non-dominated solutions obtained for 2.750 instance.}}
\end{figure}
\begin{figure}[H]
\begin{tabular}{cc}
\includegraphics[scale=0.35]{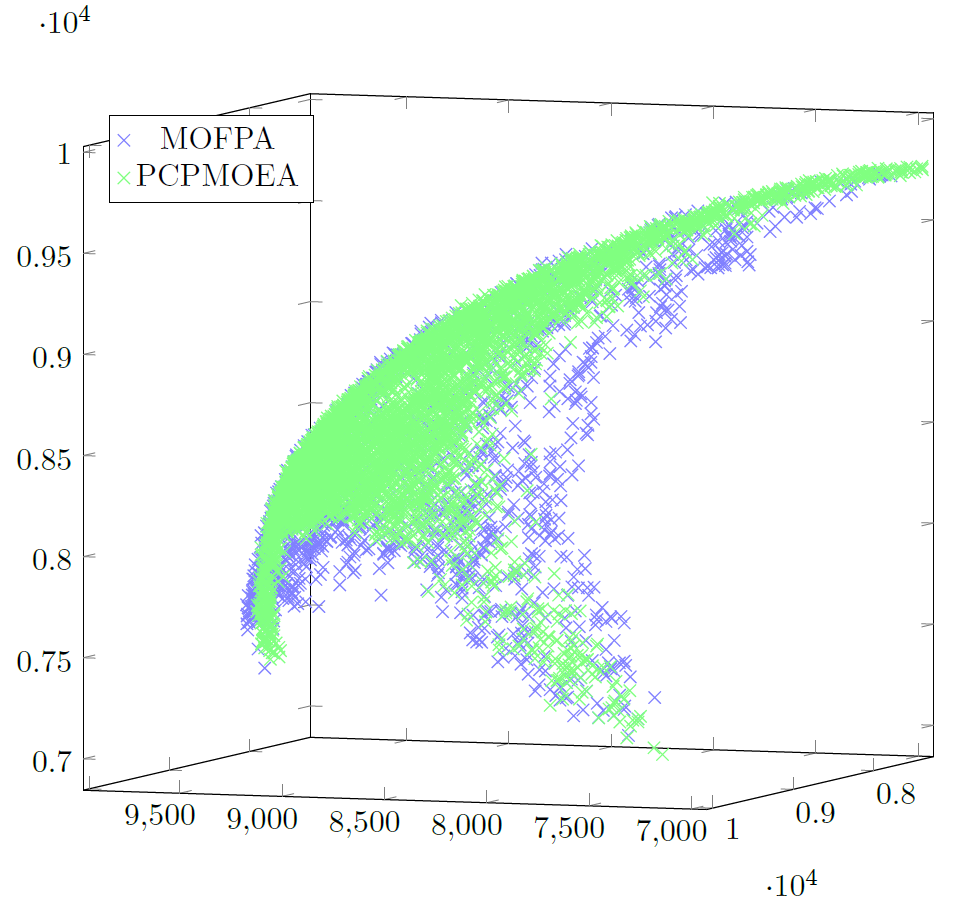} & \includegraphics[scale=0.35]{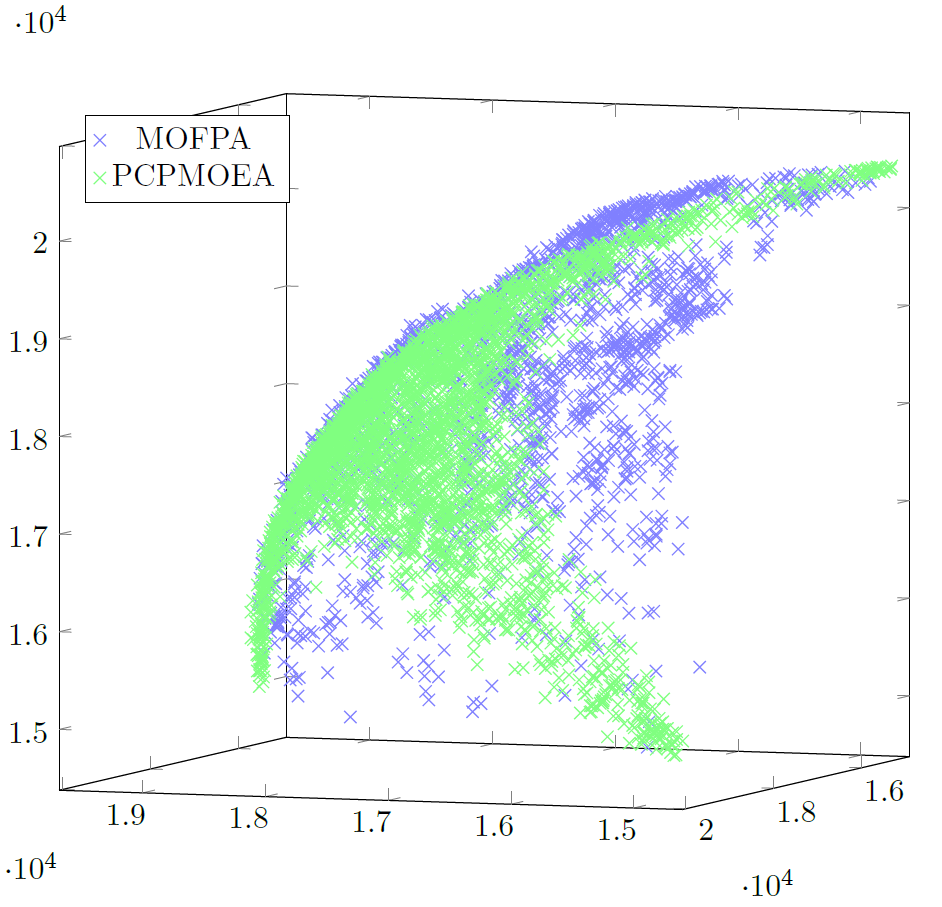}\vspace{-3mm}
\end{tabular}
\caption{{\footnotesize Illustion of Non-dominated solutions obtained for 3.250 and 3.500 instances.}}
\end{figure}

\begin{figure}[H]
\includegraphics[scale=0.37]{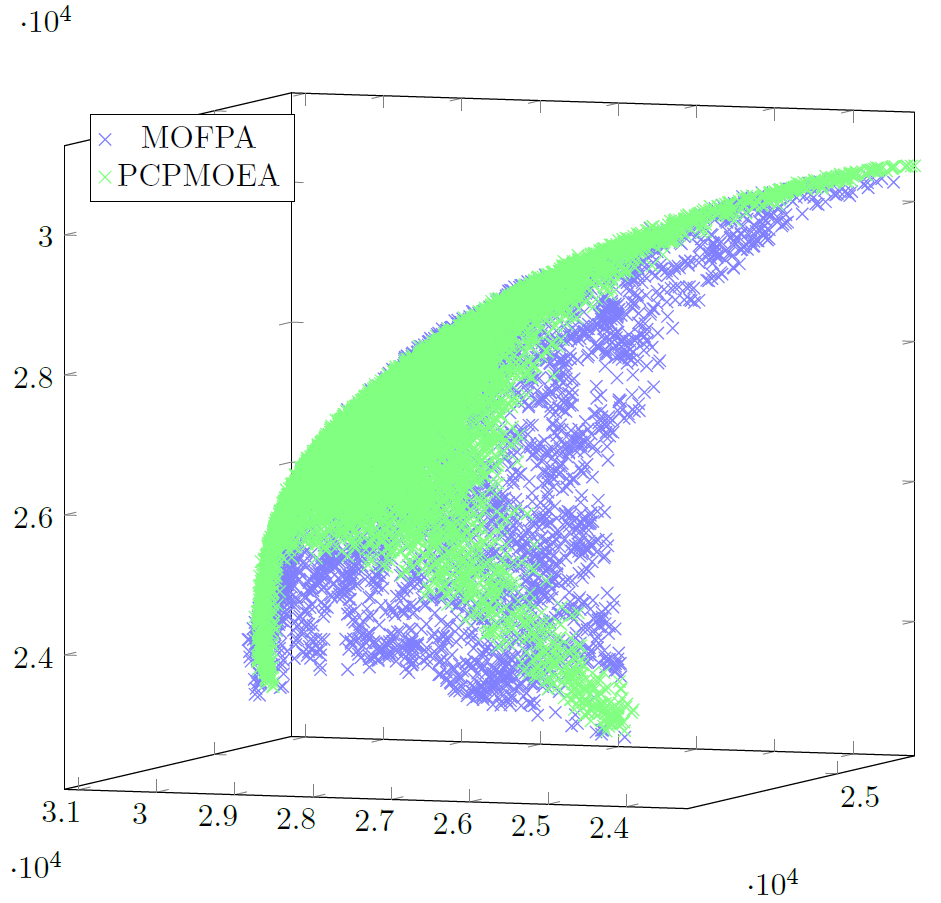}\vspace{-1mm}
\caption{{\footnotesize Illustion of Non-dominated solutions obtained for 3.750 instance.}}
\end{figure}

\section{Conclusion}

In this paper we presented a parallel multiobjective evolutionary algorithm applied to the multidimensional multiobjective Knapsack Problem (MOMKP), called Parallel Criterion-based Partitioning multiobjective evolutionary algorithm (PCPMOEA). The suggested algorithm is designed in a master/worker paradigm, handling multiple MOEAs with criterion-based selection operator. An experimental study has been carried out, where we compared the suggested algorithm, using five metrics, against four algorithms that are reputed in literature: NSGAII, SPEA2, MOEA/D, and MOFPA. Furthermore, the testes were achieved using six instances of the well-known multiobjective multidimensional Knapsack Problem, with two and three objective functions and between 250 and 500 items. The suggested algorithm has shown conclusive results regarding the standard goals of heuristic multiobjective optimization (convergence and diversity), which encourages to attempt other applications for the suggested algorithm. However, the proposed algorithm has also shown some drawbacks for large instances with three objectives. Furthermore, a brief explanation were given as we seek for possible improvements, which opens several directions to investigate as perspectives on short term future work.

\end{document}